\documentclass[aoas]{imsart}

\RequirePackage{amsthm,amsmath,amsfonts,amssymb}
\RequirePackage[authoryear]{natbib}
\RequirePackage[colorlinks,citecolor=blue,urlcolor=blue]{hyperref}
\RequirePackage{graphicx}
\usepackage{algpseudocode, multirow, threeparttable, lscape, booktabs, url, caption2}
\usepackage{algorithm}
\usepackage{mathtools}
\usepackage{color}
\usepackage{tabularx}
\usepackage{enumitem}
\usepackage{bm}

\startlocaldefs
\theoremstyle{plain}
\newtheorem{axiom}{Axiom}
\newtheorem{rmk}[axiom]{Remark}
\newtheorem{thm}{Theorem}[section]

\newtheorem{prop}[thm]{Proposition}
\newtheorem{cor}[thm]{Corollary}
\theoremstyle{definition}

\newtheorem{assump}[thm]{Assumption}

\newcommand{\R}{\mathbb{R}}

\newcommand{\RNum}[1]{\uppercase\expandafter{\romannumeral #1\relax}}

\newcommand{\lat}{{x}}
\newcommand{\lon}{{y}}
\newcommand{\Mean}{{\mathbb{E}}}
\newcommand{\Var}{{\mbox{Var}}}
\newcommand{\Cov}{{\mbox{cov}}}

\newcommand{\prob}{{\mathbb{P}}}
\def\Sigmae{{\mathbb{V}}}
\def\MSE{{\operatorname{MSE}}}
\def\ATE{{\tau}}

\newcommand{\change}[1]{\textcolor{black}{#1}}

\endlocaldefs

\begin{document}

\begin{frontmatter}
\title{Spatially Randomized Designs Can Enhance Policy Evaluation}
%\title{A sample article title with some additional note\thanksref{t1}}
\runtitle{Spatially Randomized Designs Can Enhance Policy Evaluation}
%\thankstext{T1}{A sample additional note to the title.}

\begin{aug}
%%%%%%%%%%%%%%%%%%%%%%%%%%%%%%%%%%%%%%%%%%%%%%%
%% Only one address is permitted per author. %%
%% Only division, organization and e-mail is %%
%% included in the address.                  %%
%% Additional information can be included in %%
\author[A]{\fnms{Ying}~\snm{Yang}\ead[label=e1]{yangying@fudan.edu.cn}},
\author[B]{\fnms{Chengshun}~\snm{Shi}\ead[label=e2]{C.Shi7@lse.ac.uk}}
\author[C]{\fnms{Fang}~\snm{Yao}\ead[label=e3]{fyao@math.pku.edu.cn}}
\author[D]{\fnms{Shouyang}~\snm{Wang}\ead[label=e4]{sywang@amss.ac.cn}}
\and
\author[E]{\fnms{Hongtu}~\snm{Zhu}\ead[label=e5]{htzhu@email.unc.edu}}
%%%%%%%%%%%%%%%%%%%%%%%%%%%%%%%%%%%%%%%%%%%%%%
%% Addresses                                %%
%%%%%%%%%%%%%%%%%%%%%%%%%%%%%%%%%%%%%%%%%%%%%%
\address[A]{Center for Applied Mathematics, Shanghai Key Laboratory for Contemporary Applied Mathematics, Fudan University\printead[presep={,\ }]{e1}}
\address[B]{London School of Economics and Political Science\printead[presep={,\ }]{e2}}
\address[C]{School of Mathematical Sciences, Center for Statistical Science, Peking University\printead[presep={,\ }]{e3}}
\address[D]{Academy of Mathematics and Systems Science, Chinese Academy of Sciences\printead[presep={,\ }]{e4}}
\address[E]{Department of Biostatistics, Gillings School of Global Public Health,  University of North Carolina at Chapel Hill\printead[presep={,\ }]{e5}}
\end{aug}

\begin{abstract}
This article studies the benefits of using spatially randomized experimental designs which partition the experimental area into distinct, non-overlapping units with treatments assigned randomly. Such designs offer improved policy evaluation in online experiments by providing more precise policy value estimators and more effective testing algorithms than traditional global designs, which apply the same treatment across all units simultaneously. We examine both parametric and nonparametric methods for estimating and inferring policy values based on the spatially randomized designs. Our analysis includes evaluating the mean squared error of the treatment effect estimator and the statistical power of the associated tests. Additionally, we extend our findings to the dynamic setting with spatio-temporal dependencies, where treatments are allocated sequentially over time, and account for potential temporal carryover effects. Our theoretical insights are supported by comprehensive numerical experiments.
\end{abstract}


\end{frontmatter}
%%%%%%%%%%%%%%%%%%%%%%%%%%%%%%%%%%%%%%%%%%%%%%
%% Please use \tableofcontents for articles %%
%% with 50 pages and more                   %%
%%%%%%%%%%%%%%%%%%%%%%%%%%%%%%%%%%%%%%%%%%%%%%
%\tableofcontents

\section{Introduction}
\label{sec:intro}

Policy evaluation in spatially dependent experiments involves analyzing spatially referenced data to assess the impact of new products. This methodology is widely used in diverse fields, including environmental studies \citep{zigler2012estimating}, epidemiology \citep{hudgens2008toward,callaway2023}, social science \citep{Sobel06}, and technology industries \citep{zhou2020cluster}. In such experiments, the challenge often lies in the limited number of observations and the small magnitude of treatment effects. Moreover, the strong interconnections among spatio-temporal units tend to increase the variance of estimations, making it difficult to detect the weak effects. Additionally, the treatment on one spatial unit may influence the outcomes of other units, leading to a breach of the stable unit treatment value assumption \citep[SUTVA, see, e.g.,][]{imbens2015causal}, and  causing interference or spillover effects that complicate the analysis \citep{basse2024randomization}.

As an illustration, consider the applications in ride-sourcing platforms such as Uber, Lyft, and Didi Chuxing. These companies extensively utilize A/B testing to assess the efficacy of universal treatment policies, such as new order dispatch or subsidy strategies implemented across an entire city \citep{xu2018large,zhou2021graph,luo2022policy}. In these settings, the dynamic networks of call orders and available drivers represent the supply and demand within these marketplaces, exhibiting significant spatial correlations \citep{ke2018hexagon}. Moreover, budget limitations often restrict the duration of online experiments to a mere two weeks \citep{shi2022dynamic}, with the anticipated improvements from novel policies being relatively modest, typically between 0.5\% and 2\% \citep{tang2019deep}. To understand the spillover effects, consider that implementing a subsidy policy for drivers in one area might draw drivers from adjacent units, thereby influencing outcomes in those neighboring areas as well. This interconnected structure highlights the limitations of conventional experimental designs and calls for novel frameworks tailored to such spatially and temporally dependent systems.

\begin{figure}[hb]
    \centering
    \includegraphics[width=1\textwidth,height=1in]{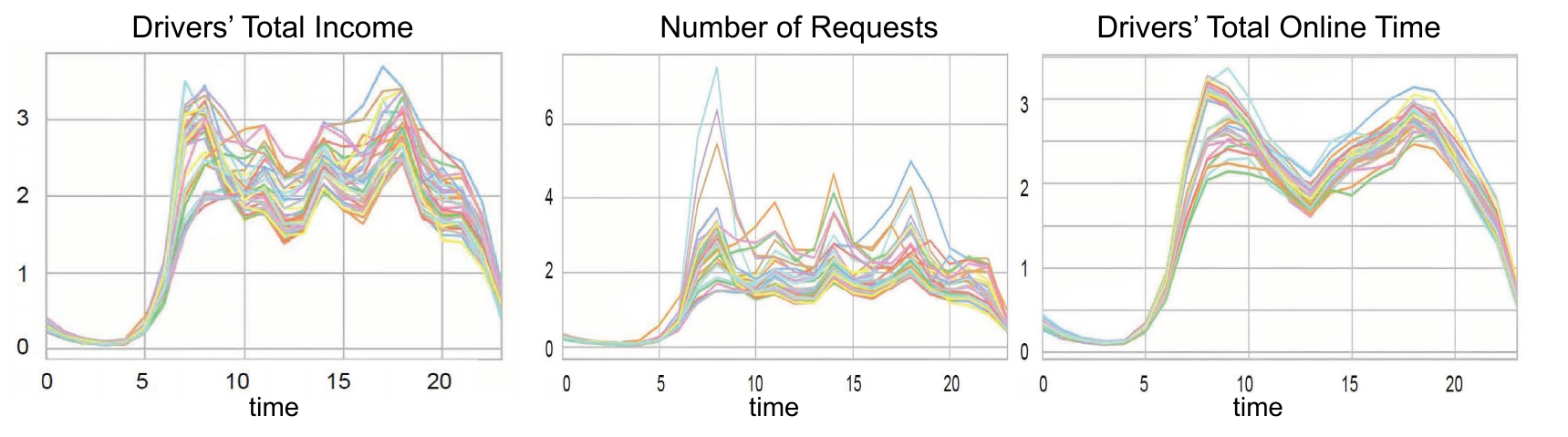}
    \vspace{-0.2in}
    \caption{ \label{fig:pattern}\small Business metrics from a city over 40 days, including drivers’ total income, the number of requests, and drivers’ total online time. Each curve represents data for a single day, with the horizontal axis corresponding to 24 hours. The values are scaled to preserve privacy.}
    \vspace{-0.15in}
\end{figure}

Experimental design plays a crucial role in accurately estimating and inferring causal effects, particularly in complex settings \citep{ding2024first,imai2013experimental,seltman2012experimental}.  In clinical trials, numerous adaptive designs have been developed to personalize treatments based on individual-level information \citep[see e.g.,][]{hu2012asymptotic,liu2022balancing,ma2024new}. However, these methods typically assume that observations are independent and identically distributed (i.i.d.), limiting their applicability in settings with spatial or temporal dependencies.
Recently, there has been growing interest in experimental designs that account for spatial dependence and interference \citep{jagadeesan2020designs,zhou2020cluster,kong2021approximate,leung2022rate}. While these works represent important advances, they generally assume independent noise across spatial units and overlook the role of temporal replicates, which is an essential feature in two-sided markets like ride-sharing platforms. In such markets, daily user and provider behaviors follow predictable patterns, such as reduced activity during late-night hours. Similar cyclical behavior is observed in other platforms, such as Airbnb and e-commerce, where interactions vary by time of day and day of week \citep{zhang2022temporal,altshuler2019modeling,liu2019exploring}. As a result, it is common to treat daily data as independent replicates. Figure~\ref{fig:pattern} illustrates these daily patterns in the ride-sharing context, reinforcing the rationale for modeling the data as temporally replicated.
Moreover, prior work largely ignores dynamic settings that capture temporal carryover effects which is an important consideration for ride-sharing companies interested in how interventions influence behavior over time. Consequently, existing literature provides limited guidance on how to design experiments in the presence of both spatial dependence and temporal replication, leaving a gap in practical methodologies for modern two-sided marketplaces.

To this end, our contributions center around the following aspects. 
Firstly, we formulate a unified framework for spatially randomized experiments that encompasses three key designs: the individual-randomized, cluster-randomized, and global designs. The global design, which applies a uniform treatment to all units, serves as a natural baseline for assessing the efficiency of spatially heterogeneous treatment allocations. To accommodate diverse data environments, we propose both parametric and semiparametric estimation procedures for evaluating treatment effects, each reflecting different levels of structural assumptions and data availability. The framework is further extended to dynamic environments characterized by spatio-temporal dependencies, enabling the study of sequential interventions and temporal spillover effects. 
Secondly, we establish quantitative results that characterize the mean squared error (MSE) and testing power of the proposed estimators, clarifying how design efficiency depends on key structural features such as spatial correlation, interference range, and cluster size. The analysis also provides insight into the relative performance of individual- and cluster-randomized designs and yields guidance on selecting the optimal cluster size that balances within- and between-cluster dependencies.. 
Thirdly, our framework directly connects to modern A/B testing and policy evaluation practices in large-scale online platforms such as ride-hailing, e-commerce, and digital advertising. In these applications, spatial or network interference naturally arises when localized interventions such as subsidies, pricing, or dispatch strategies affect outcomes in neighboring areas. The proposed spatially randomized and dynamic designs offer a principled way to account for such spillovers, enhancing both the validity and efficiency of empirical policy evaluations.
% In the realm of clinical trials, numerous adaptive designs have been developed, allowing treatments to be tailored to participants based on specific criteria and personal information \citep[see e.g.,][]{taves1974minimization,pocock1975sequential,kalish1985treatment,rosenberger2008handling,hu2012asymptotic}. These designs have proven to enhance the precision of average treatment effect estimators in traditional settings where observations are independent and identically distributed (i.i.d.), devoid of spatial or temporal interdependencies. Yet, the application and efficacy of these designs in experiments with spatial dependencies remain largely unexplored.

\subsection{Related work}\label{sec:related work}

In this subsection, we review existing works that are related to our work as follows.

\textit{\textbf{Off-policy evaluation (OPE)}}. 
 The most prevalent approach for inferring treatment effects in the aforementioned contexts is OPE within the reinforcement learning framework \citep[see e.g,][for a review]{uehara2022review}, which aims to evaluate the impact of a target policy offline using a pre-collected historical dataset generated by a different behavior policy.   
In the context of finite-horizon settings characterized by a limited number of decision points, augmented inverse propensity score weighted estimators have been introduced \citep{jiang2016doubly,luedtke2016statistical,thomas2016data}. More recent advancements have extended these methodologies to efficiently handle evaluations over extended or infinite time horizons \citep{shi2021deeply,wang2023projected,kallus2022efficiently}. 
However, these studies have not explored policy evaluation in the context of spatial interference or the role of experimental design in such settings.

\textit{\textbf{Spatial causal inference with interference}}.
Research on spatial interference has primarily branched into two prominent types of methods. The first type is the partial interference, segmenting individuals into clusters with the interference effects contained within each respective cluster \citep{Liu16,Zigler21,huber2021framework}. The second type is the local or network interference, where the interference effects are confined to the local network of each unit  \cite{perez-heydrich_assessing_2014,puelz2022graph}. 
Recent studies have proposed more complex interference structures to accommodate specific application problems \citep{Aronow17,Tchetgen21,Larsen22}. However, these works did not consider experimental designs, which is the focus of this paper.

% \textit{\textbf{Experimental designs with %spatial/network 
% interference.}}
% There is a growing interest  in experimental designs that accommodate spatial or network interference \citep{jagadeesan2020designs,zhou2020cluster,kong2021approximate,leung2022rate}. 
% \change{Their setting did not take into account replicas and assumed that the noise between spatial units is independent, nor did they consider dynamic setting, so it is not suitable for the scenario of the ride hailing market; We have considered the correlation between spatial noise and dynamic settings where temporal carryover effect exists; In terms of modeling, we consider more complex models, including both parametric and non parametric models, as well as richer theoretical results.}
% There have also been a few newly developed experimental designs that are motivated by other type of interference.
% % \cite{huber2021framework} studied the identification and estimation of spillover effects via difference-in-difference method under the cluster-randomized design. 
% \cite{bajari2021multiple,johari2022experimental,xiong2023data} studied the multiple randomization design in two-sided platforms, which is not applicable to settings with a large number of spatial units. 
% \cite{Wager21} proposed the equilibrium design to optimize supply-side payments, which is different from the target of the policy evaluation. 

\subsection{Outline of the paper}

The rest of the paper is organized as follows. In Section \ref{sec:single}, we present the problem formulation in the nondynamic setting and  utilize advanced parametric and semiparametric methods to estimate the treatment effects. We then examine the MSEs of these estimators under different spatial designs and compare their testing efficiency.  Section \ref{sec:multi-stage} extends our analysis to the dynamic case. The numerical simulations and real data applications are displayed in Section \ref{sec:expr}, which further verifies the usefulness of our method. Technical proofs are collected in the supplementary materials.

\section{Nondynamic Setting}\label{sec:single}
We start our exploration in the nondynamic setting, and the examination of the dynamic setting will be presented in the subsequent section.

\subsection{Problem formulation}
\label{subsec:DCF}

Consider a city divided into $R$ non-overlapping spatial units and we are interested in evaluating the performance of a newly developed policy against a standard control. For  the $\iota$th unit, let $Y_{\iota}(1)$ and $Y_{\iota}(0)$ represent the potential daily outcomes under the new and existing policies {for the whole city}, respectively. \change{Our focus is evaluating the average treatment effect (ATE) 
\begin{eqnarray*}
 &\ATE=\sum_{\iota=1}^R\Mean\left\{Y_{\iota}(\bm{1}_R)-Y_{\iota}(\bm{0}_R)\right\},
\end{eqnarray*} 
where $Y_{\iota}(\boldsymbol{a}_R)$ denotes the potential outcome of unit $\iota$ under the global assignment vector $\boldsymbol{A}=\boldsymbol{a}_R=(a,\ldots,a)^T\in\mathbb{R}^R$ with $a \in\{0,1\}$. This definition aligns with the policy evaluation view adopted in \cite{forastiere2021identification,leung2022rate} and \cite{lu2023design}. }
We assess the policy effectiveness through the following one-sided hypothesis test:
\begin{equation}
\label{hypo: single ols}
H_0:\ATE\le 0\text{\ \ V.S.\ \ }H_1:\ATE>0.
\end{equation}   
Under the null hypothesis, the improvement of benefits brought by the new policy is relatively low compared to the implementation costs. As such, we recommend to use the standard control. 

To estimate and test the ATE, we collect data over $N$ days in the form of treatment-outcome pairs $\{(A_{i,\iota}, Y_{i,\iota}): 1\le i\le N, 1\le \iota\le R\}$ where $A_{i,\iota}\in\{0,1\}$ denotes the binary treatment assignment for unit $\iota$ on day $i$, and $Y_{i,\iota}$ is the corresponding observed outcome. We also incorporate observed covariates $\{O_{i,\iota}\}_{i,\iota}$, such as the number of active drivers in a ride-sharing platform, to improve estimation precision.  We assume the consistency assumption (CA) as follows:

\vspace{0.05in}
\noindent\textbf{$-$CA.} For any $i$ and $\iota$, $Y_{i,\iota}$ equals to the potential outcome $Y_{i,\iota}(A_{i,\iota})$. 
\vspace{0.05in}

The aim of this paper is to compare the performance of different experimental designs (in particular, two spatially randomized designs versus a global design) in terms of policy evaluation accuracy. We evaluate their relative performance by analyzing the mean squared error (MSE) of the ATE estimators and the power of hypothesis testing procedures. To facilitate comparison, we define the efficiency ratios of the spatially randomized designs relative to the global design. We assume that observations across days are independent, i.e., daily outcomes are not influenced by treatment assignments from previous days.

\vspace{0.05in}
\textbf{Experimental designs}.
In the \textit{global design}, a uniform policy is applied across all $R$ units for each day, i.e. 
$
A_{i,1} = A_{i,2} = \cdots = A_{i,R}
$
for any given day $i$. 
In contrast, spatially randomized designs allow these $A_{i,\iota}$s to be different at each time. We study two specific types of designs. 
The first one is the \textit{individual-randomized design} that allocates treatments to each unit independently with a non-zero probability $p_{\iota}$ of receiving treatment 1.
The other one is the \textit{cluster-randomized design} which organizes units into $m$ non-overlapping clusters $\{\mathcal{C}_1,\ldots,\mathcal{C}_m\}$ based on spatial proximity and ensures uniform treatment within each cluster. Specifically, for any cluster $\mathcal{C}_j$, we have $A_{i,k_1} = A_{i,k_2} = A_j^{(j)}$ for any $k_1,k_2\in\mathcal{C}_j$, where $A_i^{(j)}$  is assigned independently with probability $p^{(j)}$.
\change{In practice, clusters are often formed by grouping spatially adjacent or strongly interacting units to reduce between-cluster interference. Typical methods include community detection or hierarchical clustering. Our theoretical analysis further suggests that that the optimal cluster size should scale with the interference range ($c^* \asymp r$), emphasizing the importance of accurately capturing the interference structure. % When this structure is unknown, it can be inferred from data using network-based detection approaches such as \citet{yuan2021causal} and \citet{zhang2024spatial}.
}
Notably, the individual-randomized design is a special case of the cluster-randomized design with $m = R$ and $\mathcal{C}_j = \{j\}$ for every $j$. 

\vspace{0.05in}
\textbf{Modeling interferences}.
In spatial settings, outcomes may be affected not only by a unit's own treatment but also by other units' treatments, a phenomenon known as interference. We model such effects via the function: $f_{\theta_\iota}\left( \{A_{i,j}\}_{j\in\mathcal{N}_\iota}\right)$ where $\mathcal{N}_\iota$ represents the interference neighbor set of unit $\iota$. That is, unit $j\in\mathcal{N}_\iota$ if the treatment on unit $j$ affects the outcome of unit $\iota$. 
These neighborhoods may be heterogeneous across units. For identifying such structures, see \cite{yuan2021causal} and \cite{zhang2024spatial}.
While various forms of $f_{\theta_\iota}(\cdot)$ can be considered such as identity functions, linear combinations, or thresholded rules, we focus on the mean-field approximation:  $\overline{A}_{i,\iota} = n_\iota^{-1} \sum_{k \in \mathcal{N}\iota} A_{ik}$, where $n_{\iota} = |\mathcal{N}_\iota|$. This formulation effectively summarizes the collective influence of neighboring treatments and has been widely adopted in the literature \citep{yang2018mean, luo2022policy, hu2022average, shi2023multiagent}. Our results can be extended to other functional forms as discussed in related work.

\subsection{Parametric and semiparametric learning}
\label{sec: sp para}
We now introduce  the parametric and semiparametric learning methods in the nondynamic setting.

\vspace{0.05in}
 \textbf{Parametric learning}. We introduce the following parametric outcome model, 
\begin{equation}
\label{model: sp para s}
Y_{i,\iota}=\alpha_\iota+O_{i,\iota}^\top\beta_\iota+\gamma_\iota A_{i,\iota}+\theta_\iota\overline{A}_{i,\iota}+e_{i,\iota},
\end{equation}
where $\alpha_{\iota}, \gamma_{\iota} \in \mathbb{R}$ and $\beta_{\iota} \in \mathbb{R}^d$. The errors $\{e_{i,\iota}\}$ are assumed to be zero-mean, temporally independent, and spatially correlated, and independent from observations and treatments. The function $f_{\theta_\iota}(\cdot)$ captures spatial spillover effects. This model allows a unit's outcome to depend on the treatments of several units while maintaining conditional independence from treatments in the other units, a common assumption in spatial analysis literature \citep[see e.g.,][]{Aronow17,reich2021review}. 
\change{We remark that one can also add the explicit neighborhood average state $\bar{O}_{i,\iota}$ in model \eqref{model: sp para s}, which is equivalent to redefining the covariate vector as $\tilde{O}_{i,\iota}=(O_{i,\iota}^T, \bar{O}_{i,\iota}^T)^T$. Since no restriction is imposed on the correlation structure among state variables, this reparametrization does not affect model generality or theoretical results.
}

Under the global design and the linearity assumption that $f_{\theta_{\iota}}(\{A_{i,j}\}_{j\in\mathcal{N}_\iota})=\theta_\iota\overline{A}_{i,\iota}$, different spatial units receive the same treatment and their outcomes satisfy
\begin{equation}
\label{model: sp para c}
Y_{i,\iota}=\alpha_\iota+O_{i,\iota}^\top\beta_\iota+\gamma_\iota^{g}A_{i}+e_{i\iota},
\end{equation}
where $\gamma_\iota^{g}=\gamma_\iota+\theta_\iota$. 
It is immediate to see that ATE has the closed-form expression $\ATE=\sum_{\iota=1}^R\gamma_\iota^{g}=\sum_{\iota=1}^R\left(\gamma_\iota+\theta_\iota\right)$.
Using data collected from the global design, we apply the ordinary least squares (OLS) regression to estimate $\gamma_{\iota}^g$ based on \eqref{model: sp para c} and plug-in these estimators to estimate the ATE, leading to  
\begin{eqnarray}
\label{eq:ols ate est single}
&\widehat\ATE^g=\sum_{\iota=1}^R\widehat{\gamma}_\iota^{g}=\sum_{\iota=1}^Ru_3^\top \left\{\sum_{i=1}^NZ_{i,\iota}^{g}\left(Z_{i,\iota}^{g}\right)^\top\right\}^{-1}\left(\sum_{i=1}^NZ_{i,\iota}^{g}Y_{i,\iota}\right),
\end{eqnarray}
where $u_3 =(0,0^T,1)^\top$ and $Z_{i,\iota}^{g}=(1,O_{i,\iota}^\top,A_{i})^\top$.

In the individual-randomized design, we adjust $\gamma_{\iota}^g$ by splitting it into $\gamma_{\iota}+\theta_{\iota}$ and apply OLS to estimate these parameters from the model. The estimators obtained are denoted as $\widehat{\gamma}_{\iota}^i$ and $\widehat{\theta}_{\iota}^i$.
For the cluster-randomized design, we define $\mathcal{C}_j^0$ as the set of ``interior'' units and $\partial \mathcal{C}_j$ as the ``boundary'' units of $\mathcal{C}_j$. A unit is in $\mathcal{C}_j^0$ if its interference neighbors are all within the same cluster $\mathcal{C}_j$. In contrast, a unit in boundary $\partial \mathcal{C}_j$ has at least one interference neighbor outside of $\mathcal{C}_j$. For these boundary units, we estimate the regression coefficients $\widehat{\gamma}_{\iota}^c$ and $\widehat{\theta}_{\iota}^c$ using OLS, similar to the individual-randomized approach. However, in the interior units, $\gamma_{\iota}$ and $\theta_{\iota}$ cannot be identified from each other, leading us to estimate their combined effect using OLS, similar to the global design.
This procedure yields the following plug-in estimators: 
\begin{eqnarray}
\label{eq:ols ate est singles}
&\widehat\ATE^{i}=\sum_{\iota=1}^Ru_{34}^\top \left\{\sum_{i=1}^NZ_{i,\iota}^{i}\left(Z_{i,\iota}^{i}\right)^\top\right\}^{-1}\left(\sum_{i=1}^NZ_{i,\iota}^{i}Y_{i,\iota}\right),\nonumber\\ &\widehat\ATE^{c}=\sum_{j=1}^m\sum_{\iota\in\mathcal{C}_j}^Ru_\iota^\top \left\{\sum_{i=1}^NZ_{i,\iota}^{c}\left(Z_{i,\iota}^{c}\right)^\top\right\}^{-1}\left(\sum_{i=1}^NZ_{i,\iota}^{c}Y_{i,\iota}\right),
\end{eqnarray}
where $u_{34}=(0,0,1,1)^\top$, $Z_{i,\iota}^i=(1,O_{i,\iota}^\top,A_{i,\iota},\overline{A}_{i, \iota})^\top$, $u_{\iota}^c=u_{34}\mathbb{I}\{\iota\in\partial\mathcal{C}_j\}+u_{3}\mathbb{I}\{\iota\in\mathcal{C}_j^0\}$,  and $Z_{i,\iota}^c=(1,O_{i,\iota}^\top,A_{i}^{(j)},\overline{A}_{i, \iota})^\top\mathbb{I}\{\iota\in\partial\mathcal{C}_j\}+(1,O_{i,\iota}^\top,A_{i}^{(j)})^\top\mathbb{I}\{\iota\in\mathcal{C}_j^0\}$.

\vspace{0.05in}
\textbf{Semiparametric learning}. We now introduce doubly robust (DR) estimators for the ATE, which are widely valued in semiparametric statistics for their resilience to model misspecification \citep[see e.g.,][]{tsiatis2006semiparametric,chernozhukov2018double}. For the outcome regression, we propose the following model: 
\begin{equation}
\label{model: sp np}
Y_{i,\iota}=h_\iota(A_{i,\iota},\overline{A}_{i,\iota},O_{i,\iota},\overline{O}_{i,\iota})+e_{i,\iota}, 
\end{equation}
where  the error terms $\{e_{i,\iota}\}_{i,\iota}$s are temporally independent, spatially correlated, and independent of the observed covariates $\{O_{i,\iota}\}$s and treatments $\{A_{i,\iota}\}$s. 

Model \eqref{model: sp np} enhances our analytical framework in two significant ways: firstly, it does not constrain the form of $h_{\iota}$, thereby permitting the use of nonparametric regression or machine learning techniques for estimation; secondly, it incorporates the influence of interference neighbor covariates $\overline{O}_{i,\iota}$, calculated as the average of $O_{i,\iota}$ across interference neighbors. This addition enriches the model by allowing more broader spatial contexts in the outcome regression model. With CA, the ATE can be expressed as 
\begin{eqnarray*}
&\ATE=\sum_{\iota=1}^R\Mean [h_\iota(1,1,O_{i,\iota},\overline{O}_{i,\iota})-h_\iota(0,0,O_{i,\iota},\overline{O}_{i,\iota})].
\end{eqnarray*}
\change{This definition reflects the effect of the target design
%rather than the experimental design used to collect the data, such as individual or cluster randomization
. The randomized experimental designs (global, individual, and cluster randomization) are introduced to enable consistent estimation of the regression functions
$h_\iota(a,a,\cdot,\cdot)$ for $a \in \{0,1\}$. % These designs determine the form of the observed data and the convergence rates for $h_\iota(a,a,\cdot,\cdot)$.
}

To construct the DR estimator, we define the following estimating function for  $1\le \iota\le R$, $1\le i\le N$ and $a\in \{0,1\}$, 
{\small\begin{equation}
\label{eqn:nuDR}
{\nu}_{DR}(a,\iota,i,h_{\iota},\pi_{\iota})=\frac{\mathbb{I}(A_{i,\iota}=a,\overline{A}_{i,\iota}=a)}{{\pi}_{\iota}(a|\{O_{i,j}\}_j)}[Y_{i,\iota}-{h}_{\iota}(a,a,O_{i,\iota},\overline{O}_{i,\iota})]+{h}_{\iota}(a,a,O_{i,\iota},\overline{O}_{i,\iota}),
\end{equation}}
where $\mathbb{I}(\bullet)$ denotes the indicator function and ${\pi}_{\iota}(a|\{O_{i,j}\}_j)$ denotes the propensity score $\mathbb{P}(A_{i,\iota}=a,\overline{A}_{i,\iota}=a|\{O_{i,j}\}_j)$ (explicit expressions are given in Section 1.1 of the supplementary materials). It follows from standard arguments that ${\nu}_{DR}(a,\iota,i,h_{\iota},\pi_{\iota})$ is an unbiased estimator for $\Mean Y_{i,\iota}(\{a\})$if either the outcome model $h_{\iota}$ or the propensity model $\pi_{\iota}$ is correctly specified. The resulting DR estimator of $\tau$ is: 
\begin{eqnarray*}
   & \frac{1}{N}\sum_{\iota=1}^R \sum_{i=1}^N [{\nu}_{DR}(1,\iota,i,\widehat{h}_{\iota},{\pi}_{\iota})-{\nu}_{DR}(0,\iota,i,\widehat{h}_{\iota},{\pi}_{\iota})],
\end{eqnarray*}
where $\widehat{h}_{\iota}$  
are estimators of $h_{\iota}$ obtained via suitable supervised learning techniques.

To avoid restrictive complexity assumptions on $h_\iota$ (e.g., metric entropy bounds \citep{diaz2020machine}), we adopt cross-fitting via sample-splitting \citep{chernozhukov2018double}. Specifically, the data are partitioned into folds; each fold is used alternately for model training and ATE estimation. The final estimator is obtained by averaging across folds, improving both robustness and efficiency. The complete procedure is detailed in Algorithm~\ref{alg:dre}.

\vspace{-0.1in}
\begin{algorithm}[htb]
\caption{Doubly robust estimation in nondynamic settings}\label{alg:dre}
\begin{algorithmic}
\Require Data $\{(O_{i,\iota}, A_{i,\iota}, Y_{i,\iota})\}_{i,\iota}$ collected from a given design.
\State \textbf{1}: Split the data into $K$ ($K\ge2$) non-overlapped subsets, each with an equal size. Let $\mathcal{I}_k$ denote the indices of the $k$th data subset. 
\State \textbf{2}: For $k=1,\ldots,K$, compute the estimated outcome regression functions $\{\widehat h_{\iota}^{(k)}\}_{\iota}$ 
based on the training data  $\{(O_{i,\iota}, A_{i,\iota}, Y_{i,\iota}): i\notin \mathcal{I}_k, 1\le\iota\le R \}$.
\State \textbf{3}: If the data is obtained from the individual-randomized design, construct the following ATE estimator 
\begin{eqnarray*}
    &\widehat\ATE_{DR}^i=\frac{1}{N}\sum_{k=1}^K\sum_{i\in \mathcal{I}_k}\sum_{\iota=1}^R\left\{\nu_{DR}(1,\iota,i,\widehat{h}_{\iota}^{(k)},\pi_{\iota}^i)-\nu_{DR}(0,\iota,i,\widehat{h}_{\iota}^{(k)},\pi_{\iota}^i)\right\}.
\end{eqnarray*}
If the data is obtained from the cluster-randomized design, construct the estimator $\widehat\ATE_{DR}^c$ by replacing $\pi_{\iota}^i$ with $\pi_{\iota}^c$ in the above equation. Otherwise, compute $\widehat\ATE_{DR}^g$ by replacing $\pi_{\iota}^i$ with $\pi_{\iota}^g$.
\end{algorithmic}
\end{algorithm}

% \vspace{-0.2in}
\subsection{Estimation accuracy in the nondynamic setting}

In this subsection, we analyze the MSEs of the ATE estimators under different models and compare different designs. 

\vspace{0.05in}
\textbf{Parametric estimators}. 
We begin by introducing the following assumption:

\vspace{-0.05in}
\begin{assump}\label{assump:omega}
Let $n_{\iota}^{(j)}$ denote the number of interference neighbors of  to the $\iota$th unit belonging to the $j$th cluster. Assume $\omega=\max_{1\le \iota\le R,\atop 1\le j_1,j_2\le m} n_{\iota}^{(j_1)}/n_{\iota}^{(j_2)}\mathbb{I}(n_{\iota}^{(j_2)}>0)=O(1)$.
\end{assump}
\vspace{-0.05in}

\begin{figure}
    \centering
    \includegraphics[width=1\textwidth,height=1.5in]{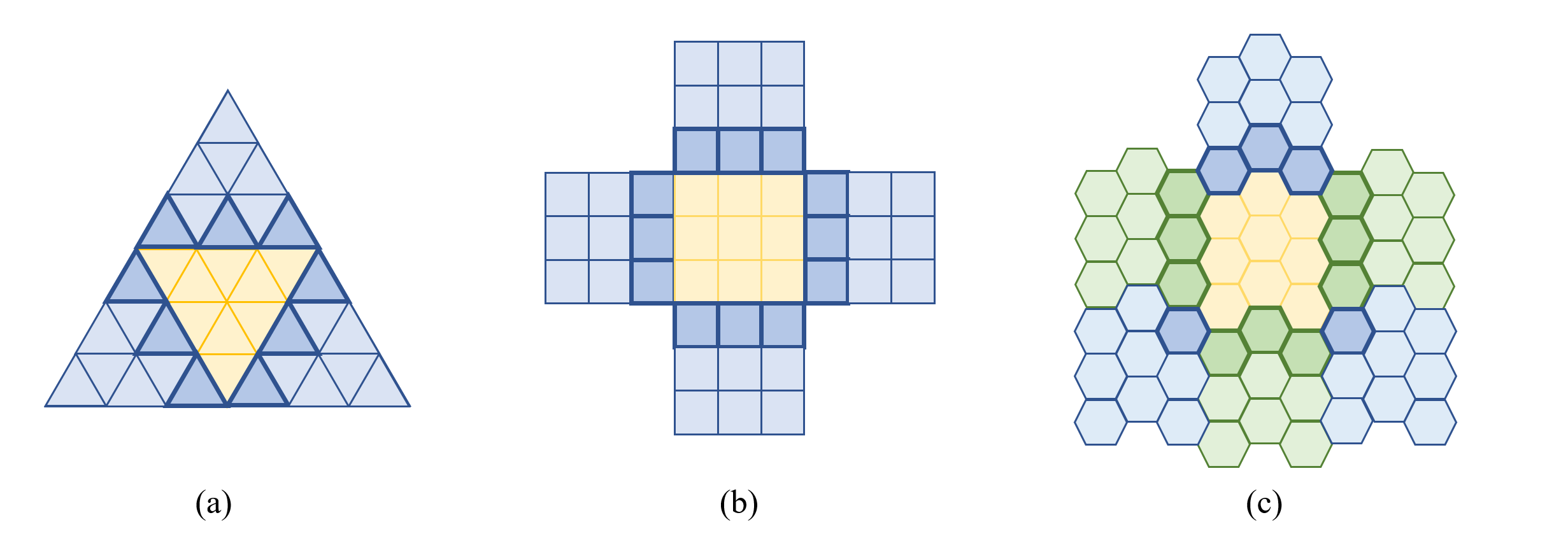}
    \vspace{-0.2in}
    \caption{ \label{fig:illustration_for_omega} \small
    Panels (a), (b), and (c) illustrate examples of clusters with interference neighbors being the adjacent units, each differing in the values of $c$ and $r$. Specifically, these panels present clusters with $c=9$ and $r=3, 4,$  and $6$, respectively. In each panel, the central cluster is emphasized in yellow, with its adjacent units outlined by bold, darker edges. Upon examination, we find that the respective values of $\omega$ and the cardinality of $|\mathcal{N}_{\mathcal{C}_j}|$ are $\omega=1, 3, 3$ and $|\mathcal{N}_{\mathcal{C}_j}|=9, 12, 14$ for panels (a), (b), and (c), respectively.
    }
    \vspace{-0.15in}
\end{figure}

Assumption~\ref{assump:omega} requires that, for any unit, the distribution of its interference neighbors across different clusters is not overly uneven, in cases where the unit and its neighbors are not all contained within the same cluster. This condition is mild and typically satisfied in practical settings--particularly when the interference range $r$ is fixed, in which case the number of neighbors per unit is uniformly bounded and the assumption holds automatically. For illustration, Figure~\ref{fig:illustration_for_omega} displays three example clusterings, where interference neighborhoods are defined based on adjacency. Let $\Sigmae=\Var(\vec{e}_i)\in\R^{R\times R}$ where $\vec{e}_i=(e_{i,1},\ldots,e_{i,R})^\top$ denotes the residual vector. We derive the MSEs of ATE estimators under arbitrary treatment probabilities $p_\iota$ and $p^{(j)}$; full details are presented in Theorem S.1 of the Supplementary Material. These results imply that MSEs are minimized when $p=p^{(j)}=p_\iota=0.5$, leading to the following theorem.

\vspace{-0.1in}
\begin{thm}\label{cor:mse single ols}
Suppose that CA holds. Set $p=p^{(j)}=p_\iota=0.5$ for all $1\le \iota\le R$, $1\le j\le m$. Let $r=\max_{\iota} n_{\iota}$ and $\nu=\sum_{\iota=1}^R \sum_{\iota'=1}^R \mathbb{V}_{\iota \iota'}/\sum_{\iota=1}^R \mathbb{V}_{\iota \iota}$, it holds that 
\begin{eqnarray}
\label{id-nondy-para}
\frac{\MSE\left(\widehat{\ATE}^i\right)}{\MSE\left(\widehat{\ATE}^g\right)}\lesssim \frac{(r+1)^2}{\nu},
\end{eqnarray}
where $a_N\lesssim b_N$ means $a_N/b_N=1+o(1)$. Further suppose that Assumption \ref{assump:omega} holds and let $\mathcal{N}_{\mathcal{C}_j}=\cup_{\iota\in\mathcal{C}_j}\mathcal{N}_\iota$. Then it holds that
\begin{equation} 
\label{cluster-nondy-para}
\frac{\MSE\left(\widehat{\ATE}^c\right)}{\MSE\left(\widehat{\ATE}^g\right)}=O\left(\frac{\sum_{j=1}^m \sum_{\iota,\iota'\in \mathcal{C}_j\cup \mathcal{N}_{\mathcal{C}_j}} \mathbb{V}_{\iota\iota'}}{\sum_{\iota,\iota'=1}^R \mathbb{V}_{\iota\iota'}}\right).
\end{equation}
\end{thm}
\vspace{-0.1in}

Theorem \ref{cor:mse single ols} characterizes how the interference structure influences the MSE of the ATE estimator in spatially randomized experiments. In particular, equation \eqref{id-nondy-para} reveals that the relative efficiency of the individual-randomized design over the   global design depends on two critical aspects: the interference range $r$ and the correlation strength $v$.  This dependency emerges because estimating spillover effects incorporates the noise from interference units into the coefficient estimation for the focal unit. When the total covariance across units scales linearly with $R$, that is,  $\nu \ge \epsilon R$ for some constant $\epsilon>0$, we obtain the bound ${\MSE(\widehat{\ATE}^i)}/{\MSE(\widehat{\ATE}^g)}=O({r^2}/{R})$. Therefore, when the interference range is small, the individual-randomized design is capable of generating a considerably more efficient estimator than the global design. 
 
In light of equation \eqref{cluster-nondy-para}, the presence of interference effects means that  $\MSE\left(\widehat{\ATE}^c\right)$ is influenced not only by the covariance matrices of residuals within the same cluster (e.g., $\{e_{\iota}\}_{\iota\in \mathcal{C}_j}$), but also by those within interference neighbors of the cluster (e.g., $\{e_{\iota}\}_{\iota\in \mathcal{N}_{\mathcal{C}_j}}$). To build intuition for the upper bound in \eqref{cluster-nondy-para}, we consider two illustrative cases:

\noindent \textbf{Case 1 (Reduction to individual randomization)}. Suppose that $m=R$ and $\mathcal{C}_j=\{j\}$ for each $1\le j\le m$, so that each unit forms its own cluster. In this case, the cluster-randomized design reduces to the individual-randomized design \citep{su2021model}. Equation~\eqref{cluster-nondy-para} then simplifies to  $O({(r+1)^2}/{\nu})$, which matches the bound in equation~\eqref{id-nondy-para}.

\noindent \textbf{Case 2 (General clustering with bounded correlation)}. Assume the pairwise covariances $\mathbb{V}_{\iota\iota'}$ are uniformly bounded above, and that the total covariance across all unit pairs, scaled by $R^2$, is bounded away from zero. Then, applying the Cauchy-Schwarz inequality, the numerator in \eqref{cluster-nondy-para} satisfies: $\sum_{j=1}^m |\mathcal{C}_j\cup \mathcal{N}_{\mathcal{C}_j}|^2\le 2\sum_{j=1}^m |\mathcal{C}_j|^2+2\sum_{j=1}^m |\mathcal{N}_{\mathcal{C}_j}|^2 $. Let $c=\max_j |\mathcal{C}_j|$ denote the largest cluster size. The first term on the right-hand side is then $O(c/R)$. Additionally, if $|\mathcal{N}_{\mathcal{C}_j}|$ is of the order of magnitude $O(r+|\mathcal{C}_j|)$, as is the case in the examples shown in Figure~\ref{fig:illustration_for_omega},  then the second term becomes $O(c/R+mr^2/R^2)$. As such, \eqref{cluster-nondy-para} is bounded by $O(c/R+mr^2/R^2)$.

This result highlights a trade-off: larger clusters increase intra-cluster correlation effects, while broader interference neighborhoods (large $r$) amplify cross-cluster spillovers. Both inflate the MSE, limiting the efficiency of the cluster-randomized design relative to the global baseline.

\vspace{0.05in}
\noindent \textbf{Semiparametric estimators}. We now examine the theoretical properties of the DR estimator for the ATE. Let $\sigma_O^2=\Var\left\{\sum_{\iota=1}^R[h_\iota(1,1,O_{i,\iota},\overline{O}_{i,\iota})-h_\iota(0,0,O_{i,\iota},\overline{O}_{i,\iota})]\right\}$ denote the variance component that is independent of the experimental design, and let $\delta_{N,\iota}^2=\max_{a\in \{0,1\} }\mathbb{E} | \widehat{h}^{(k)}_{\iota}(a,a,O_{i,\iota},\overline{O}_{i,\iota})-h_{\iota}(a,a,O_{i,\iota},\overline{O}_{i,\iota}) |^2$ represent the estimation error of the regression model for unit $\iota$.
In cluster-randomized experiments, the arrangement of a unit and its interference neighbors, whether they fall within the same cluster or are spread across multiple clusters, plays a critical role in determining the DR estimator’s performance. To capture this, we partition the residual covariance matrix according to the cluster distribution of interference neighborhoods. Let $r_c(\iota)$ denote the number of clusters that include unit $\iota$ and its interference neighbors, and define $r_c=\max_\iota r_c(\iota)$. 
Let $\mathcal{R}=\{1,2,\ldots,R\}$ as the index set of all units, and $\mathcal{R}_1=\{\iota:r_c(\iota)>2\}$. We define $$\nu_1 = \frac{\sum_\iota\sum_{\iota'}\Sigmae_{\iota\iota'}}{\sum_{\iota\in\mathcal{R}_1}\Sigmae_{\iota\iota}}, \ \nu_2=\frac{\sum_\iota\sum_{\iota'}\Sigmae_{\iota\iota'}}{\sum_{\iota\in\mathcal{R}/\mathcal{R}_1}\Sigmae_{\iota\iota}},$$ which reflect the relative covariance concentration across these two partitions.

To control the complexity introduced by interference in semiparametric models, we impose the following assumption.

\begin{assump}\label{assump:mll}
For each $\iota$, $\sum_{\iota'}\mathbb{I}(1\le m_{\iota\iota'}\le 2)=O(cr)$, where
$m_{\iota\iota'}=\sum_{k=1}^m \mathbb{I}((\mathcal{N}_{\iota}\cup \{\iota\})\cap\mathcal{C}_k \neq \emptyset, (\mathcal{N}_{\iota'}\cup \{\iota'\}) \cap \mathcal{C}_k \neq \emptyset)$
\end{assump}

This condition limits the interference overlap across clusters, ensuring that the number of unit pairs with shared or overlapping interference neighborhoods remains manageable. Assumption~\ref{assump:mll} is typically satisfied when interference is distance-based, as in the settings illustrated in Figure~\ref{fig:illustration_for_omega}. We are now ready to present the main result for the DR estimator.

{ \begin{thm}\label{cor:sp np}
Suppose that conditions in Theorem \ref{cor:mse single ols} hold, and  $R^{-1}\sum_{\iota=1}^R \delta_{N,\iota}^2\to 0$. 
% \begin{itemize}
%     \item [(i)] Let $r$ and $\nu$ be defined as they are in Theorem \ref{cor:mse single ols}.  
    Then as $N\rightarrow\infty$,
\begin{equation}
\label{eqn:firstratio}
    \frac{\MSE(\widehat\ATE_{DR}^i)-N^{-1}\sigma_O^2}{\MSE(\widehat\ATE_{DR}^g)-N^{-1}\sigma_O^2}\lesssim\frac{(r+1)^22^r}{\nu}.
\end{equation}
 When assumption \ref{assump:mll} holds, we have
 \begin{equation}
 \label{eqn:secondratio}
\frac{\MSE(\widehat{\ATE}_{DR}^c)-N^{-1}\sigma_O^2}{\MSE(\widehat{\ATE}_{DR}^g)-N^{-1}\sigma_O^2}=O\left(\frac{cr\cdot r_c\cdot 2^{r_c}}{\nu_1}+\frac{cr}{\nu_2}\right).
\end{equation}
% \end{itemize}
\end{thm}}

Theorem~\ref{cor:sp np} establishes the relative efficiency of DR estimators under individual- and cluster-randomized designs, accounting for spatial interference and estimation error in the outcome regression.
Equation~\eqref{eqn:firstratio} shows that, under a fixed interference range $r$, the MSE ratio between the individual-randomized and global designs is of order $O(R^{-1})$, indicating a potential efficiency gain. However, the ratio grows exponentially with  $r$, due to the importance sampling (IS) $\mathbb{I}(\overline{A}_{i,\iota}=a)/\prob(\overline{A}_{i,\iota}=a)$ in \eqref{eqn:nuDR}. This sensitivity to the probability of neighborhood-level treatment configurations especially as $r$ increases--is referred to as the curse of spatial interference. It represents a key limitation of the individual-randomization design in high-interference regimes. In practice, this issue can be partially mitigated by imposing a lower bound on the denominator $\prob(\overline{A}_{i,\iota}=a)$, thereby avoiding excessively large weights. 

A more principled solution is to adopt cluster-level randomization, whose benefits are reflected in the bound provided by equation~\eqref{eqn:secondratio}. The strength of cluster randomization lies in the fact that all units within a cluster receive the same treatment, reducing the variability in neighborhood treatment assignments. When Assumption~\ref{assump:mll} holds, the curse of interference is largely confined to cluster boundaries, which represent a small fraction of the spatial domain compared to the cluster interiors. The two terms on the right-hand side of equation~\eqref{eqn:secondratio} correspond to contributions from two groups of units. The first term captures the MSE of $\MSE(\widehat{\tau}^c_{DR})$ aggregated over units in $\mathcal{R}_1$, to $\MSE(\widehat{\tau}^g_{DR})$. The second term reflects the ratio of $\MSE(\widehat{\tau}^c_{DR})$, aggregated across all units in $\mathcal{R}/\mathcal{R}_1$, whose interference neighborhoods span more than two clusters. Its growth is controlled by the cluster-level interference range $r_c$, and the exponential term $2^{r_c}$ reflects the complexity of modeling these interactions. The second term accounts for the remaining units in $\mathcal{R}/\mathcal{R}_1$, which are subject to less interference across cluster boundaries. Importantly, this term lacks the exponential factor $2^{r_c}$, and is therefore generally much smaller.
These two ratios can be significantly lower than the RHS of \eqref{eqn:firstratio} for two reasons: (i) $r_c$, the cluster-level interference range, can be much smaller than $r$. For example, in the scenarios illustrated in Figure \ref{fig:illustration_for_omega}, while $r$ is 3, 4, and 6, $r_c$ is only 3, 3, and 4, respectively. (ii) The second term in \eqref{eqn:secondratio} is further diminished by a scaling factor $R^{-1}\sum_{j=1}^m |\partial \mathcal{N}_{\mathcal{C}_j}|$, which reflects the relatively small proportion of boundary units compared to interior units.

\begin{rmk}
\change{The relative efficiency between individual- and cluster-randomized designs can also be inferred from the preceding results. The MSE bounds indicate that the individual-randomized design attains smaller MSE when the number of neighboring units $r$ is small, whereas the cluster-randomized design becomes more efficient as $r$ and the within-cluster dependence increase. Specifically, under the parametric model \eqref{model: sp para s} and the conditions of Theorem \ref{cor:mse single ols}, 
$$\frac{\MSE(\widehat{\ATE}^i)}{\MSE(\widehat{\ATE}^c)}=O\left(\frac{\sum_{\iota=1}^R \mathbb{V}_{\iota \iota}}{\sum_{j=1}^m \sum_{\iota,\iota'\in \mathcal{C}_j\cup \mathcal{N}_{\mathcal{C}_j}} \mathbb{V}_{\iota\iota'}}\right).$$
When all elements of $\Sigmae_{\iota\iota'}$ are positive and uniformly bounded above and below, this simplifies to
$$\frac{\MSE(\widehat{\ATE}^i)}{\MSE(\widehat{\ATE}^c)}=O\left(\frac{r^2}{c+mr^2/R}\right).$$
Since $m \approx R / c$, it follows that $\widehat{\ATE}^c$ converges faster when $r \gg c$, whereas $\widehat{\ATE}^i$ is more efficient when $r \ll c$.
In the semiparametric case (Theorem \ref{cor:sp np}), noting that $1/\nu_1 + 1/\nu_2 = 1/\nu$, we similarly find that the cluster-randomized design is preferable when $r \gg \log c$, while the individual-randomized design is favored when $r \ll \log c$.}
\end{rmk}
\begin{rmk}
\change{Another implication of the preceding theoretical results concerns the optimal choice of cluster size in cluster-randomized designs. Under the conditions of Theorem \ref{cor:mse single ols}, when all elements of $\Sigmae_{\iota\iota'}$ are positive and uniformly bounded above and below, noting that $m \approx R / c$ and minimizing the leading term $c + m r^2 / R$, we obtain that the optimal cluster size satisfies $c^* = r$. The interference range $r$ thus serves as a key design parameter, as it characterizes the spatial extent of spillover effects and directly determines the appropriate granularity of randomization. In the semiparametric case of Theorem \ref{cor:sp np}, when the within-cluster interference radius $r_c$ remains constant, smaller clusters generally yield lower MSE. However, if $c$ becomes too small such that the cluster radius is much smaller than the interference radius (i.e., $c \ll r$), the assumption of constant $r_c$ no longer holds, and cross-cluster interference increases. This transition suggests that the optimal cluster size scales proportionally with the interference range, that is, $c^*\asymp r$, where $\asymp$ denotes asymptotic equivalence up to a constant factor.}
\end{rmk}

\subsection{Testing power in the nondynamic setting}\label{sec:test nondynamic}
In this section, we describe how to construct Wald-type test statistics for ATE inference and analyze how different experimental designs influence the power of these tests. We begin with the asymptotic normality of the estimators.

\vspace{-0.05in}
\begin{thm}\label{thm:AN}
When the conditions of Theorems \ref{cor:mse single ols} and \ref{cor:sp np}  hold,  
\begin{itemize}
\item[(2.5.1)]  there exist constants $v^g$, $v^i$ and $v^c$ 
such that 
as $N\rightarrow\infty$, 
$$\sqrt{N(v^g)^{-1}}\left(\widehat{\ATE}^g-{\ATE}\right),\ \sqrt{N(v^i)^{-1}} \left(\widehat{\ATE}^i-{\ATE}\right),\ \sqrt{N(v^c)^{-1}}\left(\widehat{\ATE}^c-{\ATE}\right)$$ 
are asymptotically standard Gaussian distributed.

\item[(2.5.2)]  there exist constants $v_{DR}^g$, $v_{DR}^i$ and $v_{DR}^c$ such that 
as $N\rightarrow\infty$, 
$$\sqrt{N(v_{DR}^g)^{-1}}\Big(\widehat{\ATE}_{DR}^{g}-\ATE\Big),\ \sqrt{N(v_{DR}^i)^{-1}}\Big(\widehat{\ATE}_{DR}^i-\ATE\Big),\ \sqrt{N(v_{DR}^c)^{-1}}\Big(\widehat{\ATE}_{DR}^c-\ATE\Big)$$ 
are asymptotically standard Gaussian distributed.
\end{itemize}
\end{thm}

Conclusion (2.5.1) follows from the classical asymptotic normality of the OLS estimator. Conclusion (2.5.2) extends the doubly robust asymptotic theory in \citet{chernozhukov2018double} to accommodate spatial interference. We now turn to the estimation of the variance terms. For the parametric estimators, standard sandwich estimators can be constructed, denoted by $\widehat{\Var}(\widehat\ATE^g), \widehat{\Var}(\widehat\ATE^i)$ and $\widehat{\Var}(\widehat\ATE^c)$. For the doubly robust estimators, the residual covariance $\Sigmae_{\iota\iota'}$ can be estimated in a similar fashion. Let $\Delta_{i\iota}=[\widehat h_\iota(1,1,O_{i,\iota},\overline{O}_{i,\iota})-\widehat h_\iota(0,0,O_{i,\iota},\overline{O}_{i,\iota})]$ and define $\overline\Delta_{i\iota}=N^{-1}\sum_{i=1}^N\Delta_{i\iota}$. The design-invariant variance component $\sigma_O^2$ can be estimated via $\widehat{\sigma}_O^2=N^{-1}\sum_{i=1}^N\{\Delta_{i\iota}-\overline{\Delta}_{i\iota}\}^2$. Alternatively, bootstrap procedures can be employed to obtain consistent variance estimates. 

To analyze the testing performance under different designs, we consider the following local alternative hypothesis:
\begin{equation}
\label{local hypo}
H_0:\ATE=0\text{\ \ V.S.\ \ }H_1:\ATE=h/\sqrt{N}.
\end{equation}
Given a consistent ATE estimator  $\widehat\ATE$ and a consistent variance estimate $\widehat\Var(\widehat\ATE)$, the standard Wald statistic is defined as $\widehat{T}={\widehat\ATE}/{\sqrt{\widehat\Var(\widehat\ATE)}}$. Let $c_\delta$ be the $(1-\delta)$th quantile of the standard Gaussian distribution and $\Phi(\cdot)$ be the cumulative density function. The null hypothesis in \eqref{local hypo} is rejected when $\widehat{T}\ge c_\delta$. Under $H_0$, Theorem~\ref{thm:AN} implies that $P\left(T\ge c_\delta\right)\rightarrow\delta$. Under $H_1$, we have the following conclusion.

\begin{cor}
\label{cor:power comp}
Suppose that the conditions of Theorems \ref{cor:mse single ols} and \ref{cor:sp np}  hold,  $\mathbb{V}_{\iota\iota'}$ is bounded from above and $\sum_{\iota\iota'}\mathbb{V}_{\iota\iota'}/R^2$ is bounded from zero. Under the alternative hypothesis in \eqref{local hypo},  as $N\rightarrow\infty$, 
\begin{itemize}
  
\item[(2.6.1)] the sufficient and necessary conditions for test consistency are:
{\small\begin{eqnarray*}
    &&P\left(\widehat{T}^g\ge c_\delta \right)\rightarrow 1\Longleftrightarrow h\gg c_\delta R,\\
    &&P\left(\widehat{T}^i\ge c_\delta \right)\rightarrow 1\Longleftrightarrow h\gg c_\delta r\sqrt{R},\\ 
    &&P\left(\widehat{T}^c\ge c_\delta \right)\rightarrow 1\Longleftrightarrow h\gg c_\delta(\sqrt{m}r+\sqrt{cR}).
\end{eqnarray*}}

\item[(2.6.2)] for the DR estimators, the corresponding conditions are:
{\small\begin{eqnarray*}
    &&P\left(\widehat{T}_{DR}^g\ge c_\delta \right)\rightarrow 1\Longleftrightarrow h\gg c_\delta\sqrt{\sigma_O^2+R^2},\\
    &&P\left(\widehat{T}_{DR}^i\ge c_\delta \right)\rightarrow 1\Longleftrightarrow h\gg c_\delta\sqrt{\sigma_O^2+r2^rR},\\ 
    &&P\left(\widehat{T}_{DR}^c\ge c_\delta \right)\rightarrow 1\Longleftrightarrow h\gg c_\delta\{\sigma_O^2+|\mathcal{R}_1|crr_c2^{r_c}+(R-|\mathcal{R}_1|)cr\}^{1/2}.
\end{eqnarray*}}
\end{itemize}
\end{cor}

These results highlight the comparative efficiency of spatially randomized designs. From (2.6.1), if the interference range $r$ is fixed and the cluster size $c\ll R$, then  both $\{r\sqrt{R}, \sqrt{m}r+\sqrt{cR}\}$ grow at a slower rate than $R$, implying that the spatial designs achieve consistency at smaller signal strengths. Similarly, in (2.6.2), when $r$ and $r_c$ are constants, $c\ll R$ and the number of boundary regions $|\mathcal{R}_1|\ll R$ (e.g., $c=O(R^\delta)$ for some $0<\delta<1$), we have 
$r2^rR\ll R^2, \text{ and } |\mathcal{R}_1|crr_c2^{r_c}+(R-|\mathcal{R}_1|)cr\ll R^2$.
Under these mild conditions, the power of tests under spatially randomized designs converges to 1 significantly faster than under the global design. This implies that, for a fixed signal strength, spatial randomization offers greater testing efficiency, particularly in large-scale experiments.

\section{Dynamic Setting}
\label{sec:multi-stage}

In this section, we consider a dynamic environment where, within each day, each spatial unit receives a sequence of $M$ (possibly different) treatments over time. The observed data take the form $\left\{\left(A_{i\iota t}, O_{i\iota t}, Y_{i\iota t}\right)\right\}_{1 \leq i \leq N,1 \leq \iota \leq R,1\le t\le M}$ where $i$ indexes the $i$th day, $\iota$ indexes the $\iota$the unit, $t$ indexes the $t$th time interval,  corresponding to i.i.d. copies of $\{\left(A_{\iota t} ,O_{\iota t} , Y_{\iota t} \right)\}_{1 \leq \iota \leq R,1\le t\le M}$. To simplify the presentation, we have chosen not to employ the potential outcome framework to formulate the causal estimand. Readers who are interested in exploring the potential outcomes approach may refer to the works by \citet{luckett2019estimating} and \citet{shi2022dynamic}. The ATE in this context is given by:
\begin{eqnarray*}
     &\ATE=\sum_{t=1}^M\sum_{\iota=1}^R [\Mean_1 (Y_{i\iota t})- \Mean_0 (Y_{i\iota t})],
\end{eqnarray*} 
where $\Mean_1$ and $\Mean_0$ denote the expectation assuming that treatments are assigned according to two non-dynamic policies: consistently set to 1 and consistently set to 0, respectively. 
In settings with multiple decision stages, past interventions may influence current variables, thereby indirectly impacting current outcomes. This carryover effects adds complexity to both the estimation process and theoretical analysis.

\subsection{Parametric and semiparametric modeling}

For parametric modeling, in addition to the outcome regression model, we posit another first-order autoregression model for the observational variable, yielding the following models: %, embodied by the following model:
\begin{equation}
\label{model:multi ols prop}
\begin{split}
&Y_{i\iota t}=\alpha_{\iota t} +O_{i\iota t}^\top\beta_{\iota t} +\gamma_{\iota t}  A_{i\iota t}+\theta_{\iota t}  \overline{A}_{i{\iota t} }+e_{i\iota t},\\
&O_{i\iota ,t+1}=\Lambda_{\iota t} + B_{\iota t} O_{i\iota t}+\Gamma_{\iota t}  A_{i\iota t}+\Theta_{\iota t}  \overline{A}_{i{\iota t} }+E_{i\iota t},
\end{split}
\end{equation}
where $\Lambda_{\iota t},\Gamma_{\iota t},\Theta_{\iota t}\in\mathbb{R}^{d}$, $B_{\iota t}\in\mathbb{R}^{d\times d}$, $\overline{A}_{i{\iota t} }=n_\iota^{-1} \sum_{k \in \mathcal{N}\iota} A_{ikt}$ and $\{e_{i\iota t},E_{i\iota t}\}$ are i.i.d. copies of $\{e_{\iota t},E_{\iota t}\}$ which are independent of $\{O_{i\iota t}\}$ and satisfy $\Mean (e_{\iota t})=0$, $\Cov(e_{\iota t},e_{\iota' t'})=\Sigmae_{\iota\iota' t}^e\mathbb{I}\{t=t'\}$,  $\Mean (E_{\iota t})=0$, and $\Cov(E_{\iota t},E_{\iota' t'})=\Sigmae_{\iota\iota' t}^E\mathbb{I}\{t=t'\}$, respectively, and for all $\iota$, ${e_{\iota t}}$ and ${E_{\iota t}}$ are independent over $t$.  This assumption rules out temporal dependence in the error structure, which, if present, may induce endogeneity and complicate identification and inference. Addressing such time-dependent endogeneity is beyond the scope of this work; see \cite{luo2022policy} for possible solutions for such settings.  
\change{In model \eqref{model:multi ols prop}, one can also include neighbor-average covariate $\bar{O}_{i\iota t}$. Though iterating the autoregressive relation would then propagate dependencies to higher-order neighbors, increasing algebraic complexity, the average treatment effect remains a linear combination of ${\gamma_{\iota t},\theta_{\iota t},\Gamma_{\iota t},\Theta_{\iota t}}_{\iota,t}$ with coefficients involving products of transition matrices. Hence the relative efficiency patterns among the three designs are expected to be qualitatively similar. For notational simplicity, we only include ${O}_{i\iota t}$ in this work.}

For the global design where $A_{i\iota t}=A_{it}$ for any $\iota$,  model (\ref{model:multi ols prop}) degenerates to
\begin{equation}
\label{model:multi ols c}
\begin{split}
&Y_{i\iota t}=\alpha_{\iota t} +O_{i\iota t}^\top\beta_{\iota t} +\gamma_{\iota t} ^{g} A_{it}+e_{i\iota t}, \\ 
&O_{i\iota ,t+1}=\Lambda_{\iota t} + B_{\iota t} O_{i\iota t}+\Gamma_{\iota t} ^{g} A_{it}+E_{i\iota t},
\end{split}
\end{equation}
where $\gamma_{\iota t} ^{g}=\gamma_{\iota t} +\theta_{\iota t}$ and $\Gamma_{\iota t} ^{g}=\Gamma_{\iota t} +\Theta_{\iota t}$.
We now turn our focus to estimating $\ATE$ within the framework of models \eqref{model:multi ols prop} and \eqref{model:multi ols c}. We first present the subsequent proposition.

\begin{prop}\label{prop:ate multi-stage linear}
Define $c_{\iota t}^\top$ as $\sum_{k=t+1}^M\beta_{\iota k}^\top\left(\prod_{j=t+1}^{k-1} B_{\iota j}\right)$. Under models \eqref{model:multi ols prop} and \eqref{model:multi ols c}, the ATE can be expressed as 
\begin{eqnarray*}
    &\ATE=\sum_{\iota=1}^R\sum_{t=1}^M(\gamma_{\iota t}^{g}+ c_{\iota t}^\top\Gamma_{\iota t}^{g})=\sum_{\iota=1}^R\sum_{t=1}^M\{\gamma_{\iota t}+\theta_{\iota t}+ c_{\iota t}^\top(\Gamma_{\iota t}+\Theta_{\iota t})\}.
\end{eqnarray*}
% 
% \begin{equation}
% % 
% \label{eqn:tauM}
% \ATE=\sum_{\iota=1}^R\sum_{t=1}^M\Big(\gamma_{\iota t}^{g}+ c_{\iota t}^\top\Gamma_{\iota t}^{g}\Big)=\sum_{\iota=1}^R\sum_{t=1}^M\Big\{\gamma_{\iota t}+\theta_{\iota t}+ c_{\iota t}^\top\big(\Gamma_{\iota t}+\Theta_{\iota t}\big)\Big\}.
% \end{equation}
\end{prop}

Proposition \ref{prop:ate multi-stage linear} shows that the ATE consists of two components: the initial terms $\gamma_{\iota t}^{c}=\gamma_{\iota t}+\theta_{\iota t}$ represents the direct effect of $A_t$ on the immediate outcomes, whereas the last term $c_{\iota t}^\top\Gamma_{\iota t}^{g}=c_{\iota t}^\top(\Gamma_{\iota t}+\Theta_{\iota t})$ corresponds to the indirect effect, measuring the delayed  treatment effect on future outcomes through $O_{t+1}$. Analogous to the nondynamic setting, one can derive OLS estimators for the regression coefficients and plug-in them into the expression in Proposition \ref{prop:ate multi-stage linear} to obtain the estimated $\tau$. To save space, we do not detail the estimating procedure again.

For semiparametric modeling, we adopt the double reinforcement learning (DRL) framework of \citet{kallus2020double}, assuming a time-varying Markov decision process (MDP) \citep{puterman2014markov}. The Markov property implies that, conditional on the current state and action $(\boldsymbol{O}_t, \boldsymbol{A}_t)$, future observations and immediate outcomes are independent of the past.
\change{Here, for $a\in\{0,1\}$, let the target policy $\pi_a$ denote the deterministic policy that sets all actions to $a$. We use $\mathbb{E}_a(\cdot)$ to denote expectation with respect to the trajectory distribution induced by $\pi_a$. The Q-function is defined as  $Q_{\iota t}^a(\boldsymbol{O}_t,\boldsymbol{A}_t)=\sum_{k=t}^M \Mean_a (Y_{\iota k}|\boldsymbol{O}_t,\boldsymbol{A}_t)$, where the current action $A_{\iota t}$ is fixed by conditioning, and the expectation is taken over future transitions $k>t$ under the policy $\pi_a$.}
We also define the density ratio between target and behavior policies as
\begin{equation*}
\mu_{\iota t}^a(\boldsymbol{O}_t,\boldsymbol{A}_t) =\frac{[\prod_{\iota} \mathbb{I}(A_{\iota t}=a)]p_{a}(\boldsymbol{O}_t)}{\prob(\cap_{\iota}\{A_{\iota t}=a\}|\boldsymbol{O}_t) p_{b}(\boldsymbol{O}_t)},\quad a\in \{0,1\},
\end{equation*}
where $p_a$ and $p_b$ denote the density functions of $\boldsymbol{O}_t$ under the target policy and the behavior policy respectively.
\change{In dynamic settings, the state distribution evolves according to the transition law $p(\bm{O}_{t+1}|\bm{O}_t,\bm{A}_t)$. Because past actions affect subsequent states, the marginal state distribution under the target policy may differ from under the behavior policy. Accordingly, the ratio  $\mu_{\iota t}^{a}$ accounts for both the difference in action probabilities and the policy-dependent state distributions.}
The resulting DRL estimator is given by 
\begin{eqnarray*}
    &\hspace{-0.6in}\widetilde{\tau}_{DRL}=\frac{1}{N}\sum_{a=0}^1 (-1)^{a+1}\sum_{i=1}^N \sum_{\iota=1}^R \Big\{Q_{\iota 1}^a(\boldsymbol{O}_{i 1}, a)\\ 
    &\hspace{1.2in}+\sum_{t=1}^M \mu_{\iota t}^a(\boldsymbol{O}_{i t},\boldsymbol{A}_{i t}) [Y_{i\iota t}+Q_{\iota t+1}^a(\boldsymbol{O}_{i t+1}, a)-Q_{\iota t}^a(\boldsymbol{O}_{i t}, \boldsymbol{A}_{it})]\Big\}.
\end{eqnarray*}
The main challenge lies in estimating $Q_{\iota t}^a$ and $\mu_{\iota t}^a$, which depend on the high-dimensional input $(\boldsymbol{O}_t, \boldsymbol{A}_t) \in \mathbb{R}^{2Rd}$. Moreover, the variance of $\widetilde{\tau}_{DRL}$ may grow exponentially with the input dimension, akin to the ``curse of horizon'' \citep{liu2018breaking}. To address this, we adopt the mean-field approximation \citep{yang2018mean,shi2023multiagent}. For each unit $\iota$, define $X_{\iota t} = (O_{\iota t}, A_{\iota t}, m_{\iota}(\boldsymbol{O}_t, \boldsymbol{A}_t))$, where $m_{\iota}(\cdot)$ summarizes local neighborhood information, e.g., averages over interference neighbors. We assume: (i) $Q_{\iota t}^a$ and $\Mean(Y_{\iota t} | \boldsymbol{O}_t, \boldsymbol{A}_t)$ depend only on $X_{\iota t}$; (ii) $O_{\iota,t+1}$ and $m_{\iota}(\cdot)$ are conditionally independent of past history given $X_{\iota t}$. These assumptions are testable using modern Markov or conditional independence tests \citep{chen2012testing, zhou2023testing, zhang2011kernel}. Under these assumptions, it can be shown that $\widetilde{\tau}_{DRL}$ is unbiased to the following \citep{shi2023multiagent}
\begin{equation*}
    \widehat{\tau}_{DRL}=\frac{1}{N}\sum_{a=0}^1 (-1)^{a+1}\sum_{i=1}^N \sum_{\iota=1}^R \Big\{Q_{\iota 1}^a(X_{i \iota t}^a)+\sum_{t=1}^M \mu_{\iota t}^a(X_{i\iota t}) [Y_{i\iota t}+Q_{\iota t+1}^a(X_{i \iota t+1}^a)-Q_{\iota t}^a(X_{i\iota t})]\Big\},
\end{equation*}
where $X_{\iota t}^a = (O_{\iota t}, a, m_{\iota}(\boldsymbol{O}_t, a))$. Unlike $\widetilde{\tau}_{DRL}$, this estimator relies on low-dimensional inputs and avoids exponential variance growth.

It remains to estimate $Q_{\iota t}^a$ and $\mu_{\iota t}^a$ to compute $\widehat{\tau}_{DRL}$. 
To estimate $Q_{\iota t}^a$, we employ backward induction \citep{murphy2003optimal}. Specifically, we begin by estimating $Q_{\iota M}^a(X_{\iota M})=\Mean (Y_{\iota M}|X_{\iota M})$. Let $\widehat{Q}_{\iota M}^a$ denote the resulting estimator. We next recursively estimate $Q_{\iota t}^a(X_{\iota t})=\Mean [Y_{\iota t}+\widehat{Q}_{\iota t+1}^a(X_{\iota t+1}^a)|X_{\iota t}]$ for $t=M-1,M-2,\cdots,1$ and construct these estimators $\{\widehat{Q}_{\iota t}^a\}_t$ via nonparametric regression 
(e.g., kernel smoothing or splines).
To estimate $\mu_{\iota t}^a$, we employ minimax learning \citep{liu2018breaking,uehara2020minimax}. To save space, we relegate the detailed implementation to Section 1.2 of the supplementary materials. 
Finally, similar to the contextual bandit setting, we employ data-splitting and cross-fitting to construct $\widehat{\tau}_{DRL}$. We summarize our procedure in Algorithm \ref{alg:dynamic}.

\begin{algorithm}[htb]
\caption{Doubly reinforcement learning in dynamic settings}\label{alg:dynamic}
\begin{algorithmic}
\Require Data $\{(O_{i,\iota,t}, A_{i,\iota,t}, Y_{i,\iota,t})\}_{i,\iota,t}$ collected from a given design.%spatially randomized designs
\State \textbf{1}: Split all data trajectories into $K$ ($K\ge2$) non-overlapped subsets, each with an equal size. Let $\mathcal{D}_k$ denote the indices of days that belong to the $k$th data subset where $k=1,2,\ldots,K$.
\State \textbf{2}: For $k=1,\ldots,K$, $a\in \{0,1\}$, $1\le t\le M$ and $1\le \iota\le R$, compute $\widehat\mu_{\iota t} ^{a,(k)}$ and $\widehat{Q}_{\iota t} ^{a,(k)}$ based on data trajectories $\{1,\ldots,n\}\backslash\mathcal{D}_k$.
\State \textbf{3}: Set 
{\footnotesize\begin{eqnarray*}
&\hspace{-0.4in}\widehat{\tau}_{DRL}=\frac{1}{N}\sum_{a=0}^1 (-1)^{a+1} \sum_{k=1}^K \sum_{i\in \mathcal{D}_k}  \sum_{\iota=1}^R \Big\{\widehat Q_{\iota 1}^{a,(k)}(X_{i \iota t}^a)+\\ 
 & \hspace{1.2in}\sum_{t=1}^M \widehat \mu_{\iota t}^{a,(k)}(X_{i\iota t}) [Y_{i\iota t}+\widehat Q_{\iota t+1}^{a,(k)}(X_{i \iota t+1}^a)-\widehat Q_{\iota t}^{a,(k)}(X_{i\iota t})]\Big\}.
\end{eqnarray*}}
\end{algorithmic}

\end{algorithm}

\vspace{-0.15in}
\subsection{Estimation accuracy in the dynamic setting}
Since treatments are sequentially assigned over time, we incorporate the spatial designs with the following temporal designs: (i) The \textit{constant design} that sets  each unit  the same treatment at each day, i.e. $A_{i\iota 1}=A_{i\iota 2}=\cdots=A_{i\iota M}$ for each $1\le \iota\le R$. (ii) The \textit{independent design} in which  all $\{A_{i\iota t}\}_t$ are independent for any $t$ and $\iota$. (iii) The \textit{switchback design} that switches the treatments for each unit back and forth (the initial treatment on each day is randomly generated), i.e. $A_{i\iota 1}=1-A_{i\iota 2}=A_{i\iota 3}=1-A_{i\iota 4}=\cdots$ for any $\iota$.

We first compare different designs under the parametric model assumption \eqref{model:multi ols prop}.  Let $\widehat\ATE_{OLS}^g, \widehat\ATE_{OLS}^i$ and $\widehat\ATE_{OLS}^c$ be the estimators under the global, individual- and cluster-randomized designs, respectively. Define the composite noise term $u_{i\iota t}=c_{\iota t}^\top E_{i\iota t}+e_{i\iota t}$, where $c_{\iota t}$ is defined the same as in Proposition \ref{prop:ate multi-stage linear} and represents the delayed indirect effects. 
Let $\Sigmae_{\iota\iota't}^u$ denote $\Cov(u_{i\iota t},u_{i\iota' t})$, and define the heterogeneity index 
\begin{eqnarray*}
     &\eta=\sum_{t=1}^M\sum_{\iota=1}^R \sum_{\iota'=1}^R \mathbb{V}_{\iota \iota' t}^u/\sum_{t=1}^M\sum_{\iota=1}^R \mathbb{V}_{\iota \iota t}^u.
 \end{eqnarray*}

\begin{thm}
\label{thm: st mse para} 
Under any of the three temporal designs (constant, independent, or switchback), it holds that
\begin{eqnarray*}
\frac{N\cdot\MSE\left(\widehat{\ATE}_{OLS}^i\right)-\sigma_{OLS}^2}{N\cdot\MSE\left(\widehat{\ATE}_{OLS}^g\right)-\sigma_{OLS}^2}\lesssim \frac{(r+1)^2}{\eta}.
\end{eqnarray*}
Furthermore, if Assumption \ref{assump:omega} holds, we have
\begin{eqnarray*}
\frac{N\cdot\MSE\left(\widehat{\ATE}_{OLS}^c\right)-\sigma_{OLS}^2}{N\cdot\MSE\left(\widehat{\ATE}_{OLS}^g\right)-\sigma_{OLS}^2}=O\left(\frac{\sum_{j=1}^m\sum_{\iota,\iota'\in\mathcal{C}_j\cup\mathcal{N}_{\mathcal{C}_j}}\Sigmae_{\iota\iota'}^u}{\sum_{\iota,\iota'}\Sigmae_{\iota\iota'}^u}\right).
\end{eqnarray*}
% \end{itemize}
\end{thm}

The results in Theorem \ref{thm: st mse para} are consistent with those in Theorem \ref{cor:mse single ols}. The explicit expression of $\sigma_{OLS}^2$ is presented in Section 3.6 of the supplementary materials.  We remark that if the indirect effects are weak, e.g., $\Gamma_{\iota t}=O(R^{-1/2})$ for any $\iota$ and $t$, or if for any $t$,   $\{E_{i\iota t}\}_{\iota}$ aare spatially correlated only within interference neighborhoods, we obtain 
\begin{equation*}
\frac{\sigma_{OLS}^2}{N\cdot\MSE\left(\widehat{\ATE}_{OLS}^g\right)-\sigma_{OLS}^2}=O\left(\frac{r^2}{R}\right).
\end{equation*}
Moreover, under standard regularity conditions, such as boundedness of $\mathbb{V}_{\iota\iota' t}^u$ and nonvanishing aggregate variance, we also have 
{\begin{equation*}
\frac{\MSE(\widehat{\ATE}_{OLS}^i)}{\MSE(\widehat{\ATE}_{OLS}^g)}=O\Big(\frac{r^2}{R}\Big)\quad\mbox{and}\quad \frac{\MSE(\widehat{\ATE}_{OLS}^c)}{\MSE(\widehat{\ATE}_{OLS}^g)}=O\Big(\frac{c}{R}+\frac{mr^2}{R^2}\Big). 
\end{equation*}}

To compare the semiparametric estimators, we assume:

\begin{assump}\label{assump:multi np}
For $1\le\iota \le R$, $1\le t\le M$ and $a=0,1$, the  reward $ Y_{\iota t}$ and the density ratio $\mu_{\iota t}^a $ are bounded; the sample splitting estimators $\widehat{\mu}_{\iota t} ^{a,(k)}$ and  $\widehat{Q}_{\iota t} ^{a,(k)}$ are finite. 
\end{assump}

Denote the DRL ATE estimators under the global, individual- and cluster-randomized designs as $\widehat\ATE_{DRL}^g,\widehat\ATE_{DRL}^i$ and $\widehat\ATE_{DRL}^c$, respectively. Let $u_{\iota t} =r_{\iota t}-\Mean(r_{\iota t}|O_{\iota t}, A_{\iota t},m_\iota(\boldsymbol{O}_{\iota t},\boldsymbol{A}_{\iota t}))$ be the noise term and $C_{\iota_1,\iota_2,t}=\Cov\{u_{\iota_1,t}, u_{\iota_2,t}\}$ be its covariance. Then we have the following result in parallel to Theorem \ref{cor:sp np}.

\begin{thm}\label{thm:st np}  
Suppose that $m_\iota(\boldsymbol{O}_t,\boldsymbol{A}_t)=(n_\iota^{-1}\sum_{k\in\mathcal{N}_\iota}O_{kt}, n_\iota^{-1}\sum_{k\in\mathcal{N}_\iota}A_{kt})$. Under the  condition that  $r$ is finite, and all elements in $\{C_{\iota_1,\iota_2,t}:1\le \iota_1,\iota_2\le R, 1\le t\le M\}$ are nonnegative with a substantial proportion distinctly greater than zero, the following result holds  as $N\rightarrow\infty$:  
\begin{align*}
\frac{N\cdot\MSE(\widehat\ATE^i_{DRL})-\sigma_{DRL}^2}{N\cdot\MSE(\widehat\ATE^{g}_{DRL})-\sigma_{DRL}^2}=O\left(\frac{r2^r}{R^2}\right).
\end{align*}
where $\sigma_{DRL}^2=\Var\{\sum_{\iota=1}^R(Q_{\iota1}(\boldsymbol{1}_{RM})-Q_{\iota1}(\boldsymbol{0}_{RM}))\}$.
Furthermore, if Assumption \ref{assump:mll} holds, then
$$\frac{N\cdot\MSE(\widehat\ATE^c_{DRL})-\sigma_{DRL}^2}{N\cdot\MSE(\widehat\ATE^{g}_{DRL})-\sigma_{DRL}^2}=O\left(\frac{cr\cdot r_c\cdot 2^{r_c}}{\nu_1}+\frac{cr}{\nu_2}\right).
$$
\end{thm}

Similar to the nondynamic case, we establish the asymptotic normality of the ATE estimators in the dynamic framework when the number of stages $M$ is finite. See Theorem S.4 in the supplementary materials for details. Comparisons of hypothesis testing performance follow directly from this result. For brevity, we omit these testing results from the main text.
\change{The same efficiency patterns between individual- and cluster-randomized designs carry over to the dynamic setting, as the underlying analytical framework and MSE decomposition remain analogous. For brevity, we omit the detailed repetition here.}

\section{Numerical Experiments}
\label{sec:expr}

In this section, we conduct numerical experiments to validate the theoretical insights derived from Sections \ref{sec:single}-\ref{sec:multi-stage}. We focus on scenarios where interference neighbors are defined as adjacent neighboring units. This setup is particularly relevant in contexts such as ride-sharing markets, where the impact of a policy in one area can influence neighboring units through driver distribution, under the assumption that drivers can only move to adjacent areas within a given time frame.

\vspace{-0.05in}
\subsection{Simulation of the nondynamic setting}\label{sec:simu single}

\begin{figure}[htbp]
\centering
\includegraphics[width=1\textwidth,height=2.6in]{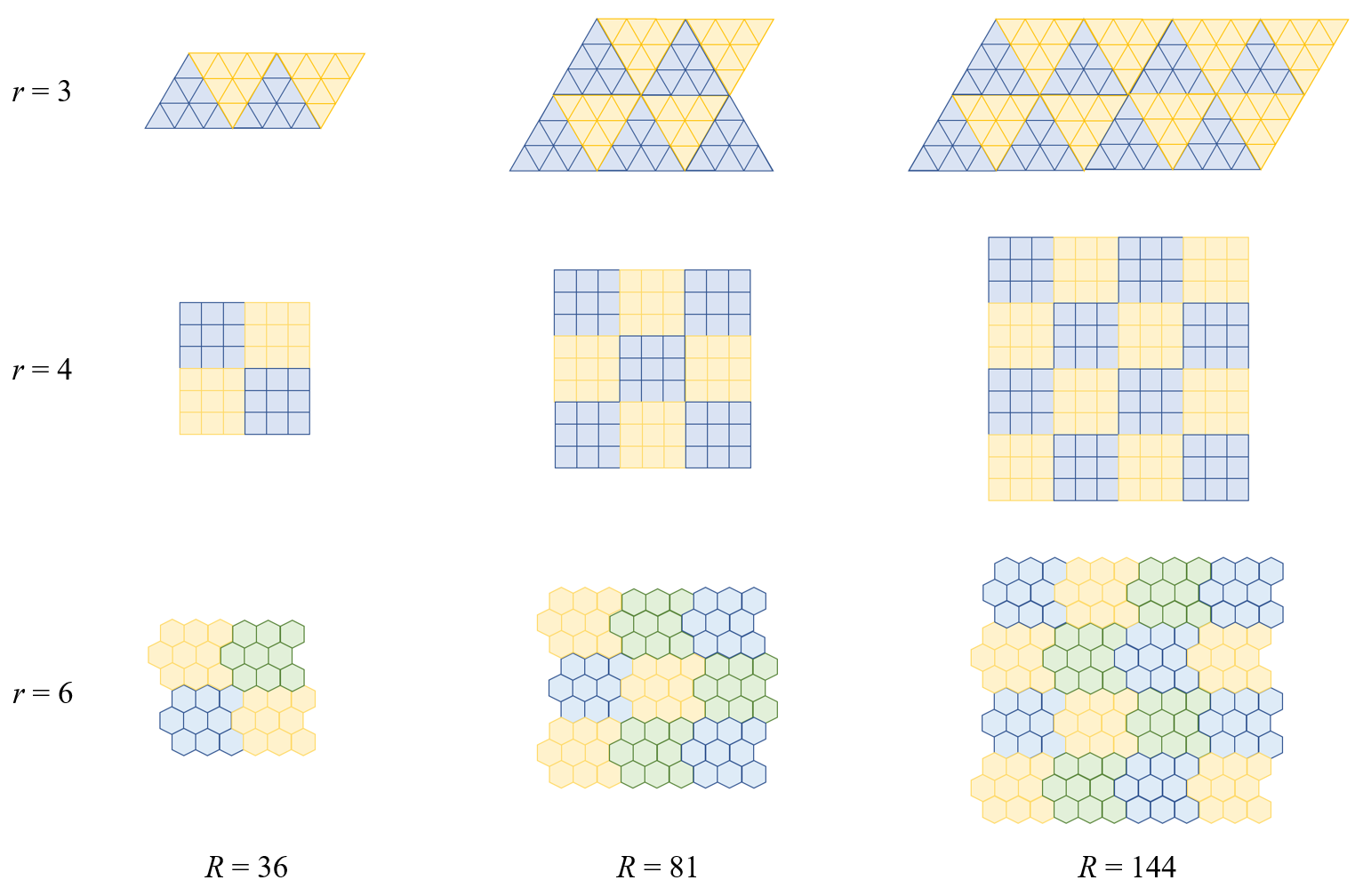}
    \vspace{-0.2in}
\caption{ \label{fig:unit_pattern}\small Simulation layout showcasing unital patterns. Each row illustrates configurations with a maximum number of neighbors per unit, specifically $r=$3, 4, and 6. The columns indicate varying total counts of units, with $R=$36, 81, and 144. Within each panel, distinct colors denote separate clusters, leading to varying cluster counts of $m=$4, 9, and 16. }
    \vspace{-0.15in}
\end{figure}

We evaluate performance in a single-stage setting with $R=36$, $81$, and $144$ units, arranged as regular polygons with $r=3$, $4$, and $6$ sides, respectively. Clusters are of fixed size $|\mathcal{C}_j|=9$, yielding $m=4$, $9$, and $16$ clusters  as depicted in Figure~\ref{fig:unit_pattern}. 

In the parametric model \eqref{model: sp para s}, covariates $O_{i,\iota}$ are i.i.d  following $N(4,1)$. Each unit $\iota$ has normalized coordinates $(\iota_x, \iota_y) \in (0,1)^2$. We set $\alpha_{\iota} = 8 + 2\{f_1(\iota_x) + g_1(\iota_y)\}$ and $\beta_{\iota} = f_2(\iota_x) + g_2(\iota_y)$, where $f_k, g_k$ are Fourier series of the form $a_0 + \sum_{k=1}^K (a_k \cos k\pi x + b_k \sin k\pi x)$ with $K=3$ and coefficients drawn from $U(0,1)$. Treatment effects follow $\gamma_{\iota} = s_y \cdot \alpha_{\iota} / \sum_{\iota} \alpha_{\iota}$ and $\theta_{\iota} = 0.6s_y \cdot \beta_{\iota} / \sum_{\iota} \beta_{\iota}$, where $s_y$ corresponds to $s\%$ improvement over the control outcome. We vary $s \in \{0, 0.25, 0.5,\dots, 2\}$.

For the semi-parametric model, the underlying model is 
\begin{equation*}
Y_{t,\iota}=5+3(O_{i,\iota}+\overline{O}_{i\iota})\sin\left\{\frac{\pi}{8}(\iota_x+\iota_y)+s\cdot A_{t,\iota}+0.5\cdot s \cdot \overline{A}_{i,\iota}\right\} + 0.5e_{i\iota},
\end{equation*}
where $O_{i,\iota} \sim N(4,1)$ truncated to $(3,5)$ and $s$ varies from $0$ to $0.015$ in increments of $0.0025$. The noise terms $e_{i\iota}$ are mean-zero with $\Cov(e_\iota, e_{\iota'}) = \rho^{d_{\iota\iota'}}$, where $d_{\iota\iota'} = \sqrt{(\lat_\iota - \lat_{\iota'})^2 + (\lon_\iota - \lon_{\iota'})^2}/2$ is the scaled spatial distance. We vary $\rho \in \{0.3, 0.6, 0.9\}$, with $N=30$ samples per configuration and 500 replications. Additional noise structures are discussed in Section 4 of the supplementary materials.

Table~\ref{tab:mse single spill} presents empirical MSE ratios, where $r_1$ and $r_2$ denote the ratios of the individual- and cluster-randomized designs to the global design, respectively, for both parametric ($\widehat\ATE$) and semiparametric ($\widehat\ATE_{DR}$) estimators. The results demonstrate that spatially randomized designs consistently reduce MSEs compared to the global design, especially as the noise correlation $\rho$ and/or the number of units $R$ increases. The individual-randomized design is generally more efficient than the cluster-randomized one when the number of neighbors $r$ is small (e.g., $r=3,4$), while the cluster-randomized design becomes preferable for larger $r$ (e.g., $r=6$). 

\begin{table}[htbp]
  \centering
  \caption{ \label{tab:mse single spill} Nondynamic setting: MSE ratio comparisons under parametric and semiparametric models across different $(r,\rho)$ settings.}
\resizebox{1\textwidth}{0.9in}{\begin{tabular}{|c|c|c|ccccccccc|}
    \hline
    Model & $R$ & Ratio & (3,0.9) & (3,0.6) & (3,0.3) & (4,0.9) & (4,0.6) & (4,0.3) & (6,0.9) & (6,0.6) & (6,0.3) \\
    \hline
    \multirow{6}{*}{Parametric}
     & \multirow{2}{*}{36} & $r_1$ & 0.360 & 0.379 & 0.421 & 0.544 & 0.582 & 0.661 & 0.982 & 1.042 & 1.171 \\
     &                     & $r_2$ & 0.522 & 0.554 & 0.593 & 0.641 & 0.690 & 0.759 & 0.674 & 0.724 & 0.795 \\
     \cline{2-12}
     & \multirow{2}{*}{64} & $r_1$ & 0.214 & 0.200 & 0.195 & 0.316 & 0.306 & 0.319 & 0.577 & 0.557 & 0.577 \\
     &                     & $r_2$ & 0.370 & 0.361 & 0.361 & 0.457 & 0.460 & 0.485 & 0.493 & 0.497 & 0.524 \\
     \cline{2-12}
     & \multirow{2}{*}{144}& $r_1$ & 0.117 & 0.130 & 0.155 & 0.187 & 0.219 & 0.259 & 0.357 & 0.413 & 0.485 \\
     &                     & $r_2$ & 0.185 & 0.196 & 0.212 & 0.246 & 0.272 & 0.090 & 0.271 & 0.298 & 0.339 \\
    \hline
    \multirow{6}{*}{Semiparametric}
     & \multirow{2}{*}{36} & $r_1$ & 0.186 & 0.202 & 0.239 & 0.352 & 0.382 & 0.459 & 0.286 & 0.316 & 0.385 \\
     &                     & $r_2$ & 0.306 & 0.331 & 0.378 & 0.438 & 0.484 & 0.576 & 0.363 & 0.399 & 0.466 \\
     \cline{2-12}
     & \multirow{2}{*}{64} & $r_1$ & 0.103 & 0.106 & 0.118 & 0.165 & 0.175 & 0.204 & 0.138 & 0.153 & 0.185 \\
     &                     & $r_2$ & 0.163 & 0.173 & 0.192 & 0.199 & 0.214 & 0.244 & 0.190 & 0.203 & 0.231 \\
     \cline{2-12}
     & \multirow{2}{*}{144}& $r_1$ & 0.062 & 0.068 & 0.079 & 0.096 & 0.107 & 0.127 & 0.080 & 0.090 & 0.107 \\
     &                     & $r_2$ & 0.089 & 0.095 & 0.107 & 0.102 & 0.111 & 0.129 & 0.102 & 0.111 & 0.129 \\
    \hline
  \end{tabular}}
\end{table}%

\begin{figure}[h]
\centering
\includegraphics[width=1\textwidth,height=2.6in]{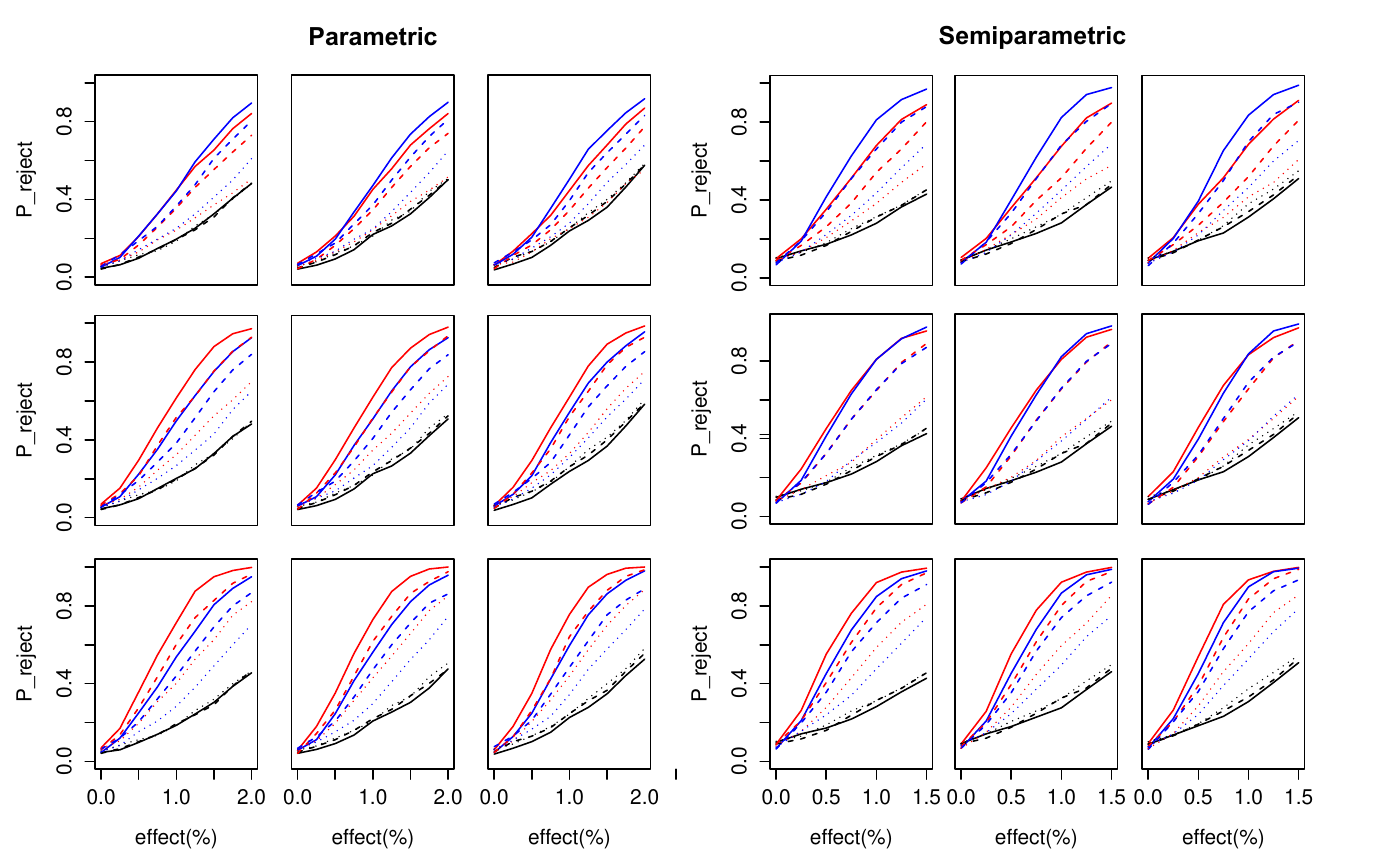}
    \vspace{-0.2in}
\caption{ \label{fig:single_mode2} \small
Rejection probabilities in the parametric and semiparametric models under a non-dynamic setup. The horizontal axis represents the treatment's relative improvement. The red lines depict the individual-randomized design,  the blue lines represent the cluster-randomized design, and the black lines denote the global randomized design. The line styles-solid, dashed, and dotted-correspond to the number of units $R=144, 81,$ and $36,$ respectively. The figure is organized into three rows and three columns of panels, representing different values of $r$ (6, 4, and 3) and $\rho$ (0.9, 0.6, and 0.3), respectively.
}
    \vspace{-0.15in}
\end{figure}

Inference results based on these settings show that both designs control type I error near the nominal 0.05 level when $s=0$. As the treatment effect increases, spatially randomized designs exhibit superior power under both models. Figure~\ref{fig:single_mode2} further confirms that power improves with larger $R$ or smaller $r$, aligning with our theoretical findings.

\subsection{Simulation of the dynamic setting}\label{sec:simu multi}

In the dynamic setting, we use the same unit configurations as in the static case. For the parametric model, initial covariates $O_{i\iota 1} \sim N(4,1)$ are i.i.d., and for $t \geq 1$, covariates $O_{i\iota t}$ follow model \eqref{model:multi ols prop}. Coefficients $\alpha_{\iota t}, \beta_{\iota t}, \Lambda_{\iota t}$ vary over time via functions $h_1(t), h_2(t), h_3(t)$:
\begin{equation*}
\alpha_{\iota t}=h_1(t)\{8+2(f_1(\iota_x)+g_1(\iota_y))\},\ \beta_{\iota t}=h_2(t)\{f_2(\iota_x)+g_2(\iota_y)\},\ \Lambda_{\iota t}=h_3(t)\{f_3(\iota_x)+g_3(\iota_y)\},
\end{equation*}
where each $h_i(t)$ is of the form $a_0 + \sum_{k=1}^K (a_k\cos(k\pi t)+b_k\sin(k\pi t))$ with $a_0,a_k,b_k \sim U(0,1)$, $K=3$. Similarly, define $B^0_{\iota t} = h_4(t){f_4(\iota_x)+g_4(\iota_y)}$, then normalize: $B_{\iota t} = 0.5(B^0_{\iota t} - \min B^0)/(\max B^0 - \min B^0) + 0.3$.
Noises $e_{i\iota t}$ are i.i.d. as in the static case, while $E_{i\iota t}$ are i.i.d. copies of $ 0.1e_\iota$. Effect coefficients are scaled by signal strengths $s_x$ and $s_y$: $\gamma_{\iota t}=s_y\cdot{\alpha_{\iota t}}/{\sum_\iota\sum_t\alpha_{\iota t}}$, $\theta_{\iota t}=0.6s_y\cdot\cdot{\beta_{\iota t}}/{\sum_\iota\sum_t\beta_{\iota t}}$,  $\Gamma_{\iota t}=s_x\cdot{\Lambda_{\iota t}}/{\sum_\iota\sum_t\Lambda_{\iota t}}$ and $\Theta_{\iota t}=0.6s_x\cdot\cdot{B_{\iota t}}/{\sum_\iota\sum_t B_{\iota t}}$,  with $s_x = s \cdot \sum_{\iota,t} \Mean(O_{\iota t}|A_{jk}=0)$ and $s_y = s\% \cdot \sum_{\iota,t} \Mean(Y_{\iota t}|A_{jk}=0)$ representing $s\%$ improvements.

For the nonparametric model, the underlying model is 
\begin{equation*}
Y_{i\iota t}=5+2\cdot O_{i\iota t}\sin\left\{\frac{\pi}{8}(\iota_x+\iota_y+t/M)+s\cdot (A_{i\iota t}+ \overline{A}_{i\iota t})\right\} + 0.5e_{i\iota t},
\end{equation*}
where $(O_{i1t},\ldots,O_{iRt})^\top$ are i.i.d. multivariate Gaussians with mean $(2+A_{i1t},\ldots,2+A_{iRt})^\top$ and covariance $\Sigma(\rho))$ having unit variances and off-diagonals $v_\iota v_\iota^T$, where $v_\iota\sim \rm{Unif}(0.75,1)$ for about a $\rho$-.

\begin{figure}[h]
\centering
\includegraphics[width=0.8\textwidth,height=2in]{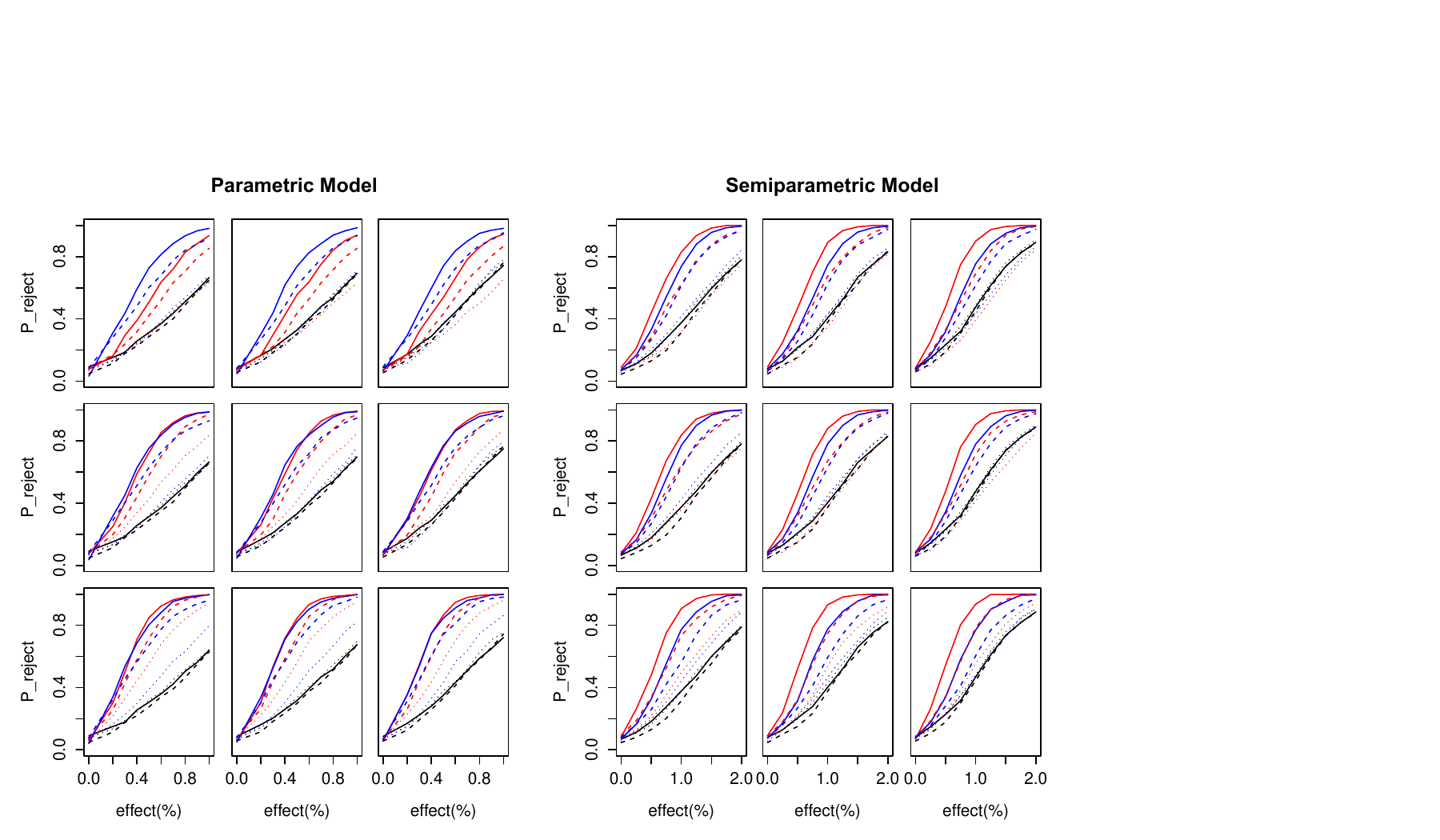}
    \vspace{-0.15in}
\caption{ \label{fig:multi} \small Rejection probability in the parametric and nonparametric models with $M=12$ of dynamic setting under different relative improvements of the new policy. The red, blue and black lines represent to the individual-, cluster- and global-randomized designs, with $R=144,81,36$ plotted in solid, dashed and dotted lines, respectively. The three rows of panels correspond to $r=6,4,3$, and the three columns correspond to $\rho=0.9,0.6,0.3$, respectively.}
    \vspace{-0.15in}
\end{figure}

We simulate $r=3,4,6$ (number of neighbors), $\rho = 0.3, 0.6, 0.9$ (noise correlation), $R = 36, 81, 144$ (units), and $M = 12$ time steps, each with $N=30$ observations. Table~\ref{tab:multi} summarizes empirical MSE ratios across settings, showing consistent trends with the static case and supporting Theorems~\ref{thm: st mse para} and~\ref{thm:st np}. Figure~\ref{fig:multi} presents power curves for various $s$ levels. In all scenarios, spatially randomized designs outperform the global design, highlighting their effectiveness in dynamic environments fraction of indices and 0 otherwise, yielding a low-rank correlation structure controlled by $\rho$.

\begin{table}[htbp]
  \centering
  \caption{ \label{tab:multi} Dynamic setting: MSE ratio comparisons under parametric and semiparametric models across different $(r,\rho)$ settings.}
\resizebox{1\textwidth}{0.9in}{\begin{tabular}{|c|c|c|c|c|c|c|c|c|c|}
    \hline
    Model & $R$ & Metric & (3, 0.9) & (3, 0.6) & (3, 0.3) & (4, 0.9) & (4, 0.6) & (4, 0.3) & (6, 0.9) \\
    \hline
    \multirow{6}{*}{Parametric}
      & \multirow{2}{*}{36} & $r_1$ & 0.315 & 0.350 & 0.410 & 0.551 & 0.603 & 0.694 & 1.067 \\
      &                     & $r_2$ & 0.495 & 0.513 & 0.543 & 0.693 & 0.738 & 0.808 & 0.726 \\
      \cline{2-10}
      & \multirow{2}{*}{64} & $r_1$ & 0.249 & 0.252 & 0.267 & 0.434 & 0.444 & 0.475 & 0.709 \\
      &                     & $r_2$ & 0.291 & 0.302 & 0.326 & 0.388 & 0.416 & 0.471 & 0.428 \\
      \cline{2-10}
      & \multirow{2}{*}{144} & $r_1$ & 0.116 & 0.126 & 0.142 & 0.196 & 0.211 & 0.238 & 0.366 \\
      &                      & $r_2$ & 0.103 & 0.114 & 0.138 & 0.158 & 0.179 & 0.224 & 0.165 \\
    \hline
    \multirow{6}{*}{Semiparametric}
      & \multirow{2}{*}{36} & $r_1$ & 0.139 & 0.135 & 0.152 & 0.181 & 0.156 & 0.181 & 0.193 \\
      &                     & $r_2$ & 0.641 & 0.713 & 0.809 & 0.666 & 0.757 & 0.871 & 0.705 \\
      \cline{2-10}
      & \multirow{2}{*}{64} & $r_1$ & 0.073 & 0.070 & 0.075 & 0.084 & 0.082 & 0.085 & 0.089 \\
      &                     & $r_2$ & 0.425 & 0.474 & 0.587 & 0.355 & 0.409 & 0.503 & 0.352 \\
      \cline{2-10}
      & \multirow{2}{*}{144} & $r_1$ & 0.033 & 0.032 & 0.033 & 0.037 & 0.035 & 0.038 & 0.043 \\
      &                      & $r_2$ & 0.224 & 0.250 & 0.296 & 0.218 & 0.245 & 0.291 & 0.221 \\
    \hline
  \end{tabular}}
\end{table}%

\subsection{Real data based simulation}

In this section, we evaluate the performance of our proposed experimental framework by simulating a real online ride-hailing service, using the same simulator  in \cite{zhou2021graph}. This simulator accurately models the dynamic interplay between demand and supply, mirroring an actual ride-hailing platform. It generates demand distributions from historical data, while initial driver supply distributions are based on early-day data, subsequently adjusted through the simulator's dynamic transitions and order dispatch algorithms. The simulation's key metrics, such as driver earnings, response rates, and idle times, deviate from real-world figures by less than 2\%.

For our analysis, we selected the largest square unit of the city, dividing it into $64$ equally sized, non-overlapping squares ($R=64$). Over $N=20$ days, segmented into $M=24$ hourly intervals, we investigated the effects of policies simulating ATE improvements of $0\%, 1.0\%, 2\%$, and $5\%$. We applied a parametric regression model for the multi-stage scenario to estimate and analyze the ATE, comparing the classical global design with an individual randomized design that incorporates temporal randomization. The outcomes of these designs are denoted by $\widehat\ATE^i_{OLS}$ and $\widehat\ATE^g_{OLS}$, respectively.

Using 300 Monte Carlo simulations to closely estimate the true ATE, we compiled the test powers and MSE ratios $\MSE(\widehat\ATE^i_{OLS})/\MSE(\widehat\ATE^g_{OLS})$ in Table \ref{tab:simulator}. Our findings indicate that both the classical global and spatially randomized designs maintain proper type I error rates. Notably, the spatially randomized design consistently exhibited lower MSE and higher test power compared to the classical global approach, irrespective of interference effects. This evidence strongly supports the superior efficiency of the spatially randomized design in real-world applications.

% \begin{table}[htbp]
% \centering
% \caption{ \label{tab:simulator}The test powers and  $\MSE(\widehat\ATE_{OLS}^i)/\MSE(\widehat\ATE_{OLS}^g)$ of the data produced by the Simulator.}
% \begin{tabular}{|c|l|r|r|r|}
% \hline
% \multicolumn{2}{|c|}{\multirow{2}[4]{*}{}} & \multicolumn{2}{c|}{ATE-test} & \multicolumn{1}{c|}{\multirow{2}[1]{*}{MSE ratio}} \\
% \cline{3-4}    \multicolumn{2}{|c|}{} & \multicolumn{1}{c|}{benchmark} & \multicolumn{1}{c|}{proposed} &  \\
% \hline
% \multirow{2}[0]{*}{0\%} & no-interference & \multirow{2}[0]{*}{0.07 } & 0.06  & 0.010 \\
% \cline{2-2}\cline{4-5}          &  interference-existing &       & 0.05  & 0.129 \\
% \hline
% \multirow{2}[0]{*}{1.0\%} & no-interference & \multirow{2}[0]{*}{0.16 } & 1.00     & 0.013 \\
% \cline{2-2}\cline{4-5}          &  interference-existing &       & 0.32  & 0.134 \\
% \hline
% \multirow{2}[0]{*}{2\%} & no-interference & \multirow{2}[0]{*}{0.33 } & 1.00     & 0.016 \\
% \cline{2-2}\cline{4-5}          &  interference-existing &       & 0.89  & 0.128 \\
% \hline
% \multirow{2}[0]{*}{5\%} & no-interference & \multirow{2}[0]{*}{0.87 } & 1.00     & 0.037 \\
% \cline{2-2}\cline{4-5}          &  interference-existing &       & 1.00     & 0.158 \\
% \hline
% \end{tabular}
% \end{table}%

\begin{table}[htbp]
\centering
\caption{\label{tab:simulator}Test powers and MSE ratios $\MSE(\widehat\ATE_{OLS}^i)/\MSE(\widehat\ATE_{OLS}^g)$ under different effect sizes and interference settings.}
\renewcommand{\arraystretch}{1.2}
\begin{tabular}{|c|cc|cc|cc|}
\hline
\multirow{2}{*}{Effect Size} 
& \multicolumn{2}{c|}{ATE-test (No Interference)} 
& \multicolumn{2}{c|}{ATE-test (With Interference)} 
& \multicolumn{2}{c|}{MSE Ratio} \\
\cline{2-7}
& Benchmark & Proposed & Benchmark & Proposed & No Interf. & Interf. \\
\hline
0\%   & 0.07 & 0.06 & 0.07 & 0.05 & 0.010 & 0.129 \\
1\%   & 0.16 & 1.00 & 0.16 & 0.32 & 0.013 & 0.134 \\
2\%   & 0.33 & 1.00 & 0.33 & 0.89 & 0.016 & 0.128 \\
5\%   & 0.87 & 1.00 & 0.87 & 1.00 & 0.037 & 0.158 \\
\hline
\end{tabular}
\end{table}

\vspace{-0.05in}
\subsection{Real data example}

\change{To illustrate our method, we analyze data from a dynamic A/B experiment conducted on the Didi Chuxing ride-sourcing platform. The experiment tested a customer subsidy policy in a metropolitan area partitioned into $R=17$ regions. Treatments were assigned using an individual-randomized switchback design over $N=24$ days (2020/02/19–2020/03/13) with 30-minute intervals ($M=48$). Each region’s adjacent regions were defined as its interference neighbors, a common assumption in spatial transport analysis. Over 90\% of ride requests were served by drivers within the same region, ensuring primarily local interactions. The number of requests was used as the state variable, and drivers’ total income was the outcome. Although the number of spatial units in this example is modest while the number of time intervals is relatively large, the setting still provides a meaningful illustration of our method in practice. The analysis demonstrates how spatially randomized designs can be implemented and evaluated in realistic dynamic environments, even when the spatial dimension is limited.}

\begin{figure}[h]
\centering
\includegraphics[width=1\textwidth]{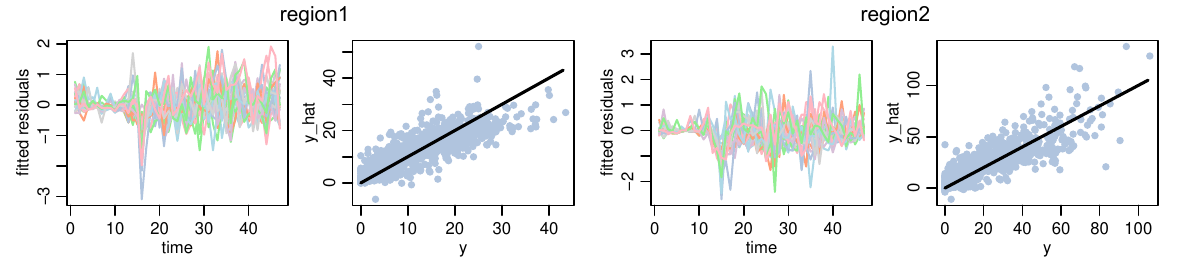}
    \vspace{-0.2in}
\caption{ \label{fig:realdata} \small Plots of the corresponding relative residuals, as well as the fitted drivers’ income against the observed values.}
    \vspace{-0.15in}
\end{figure}

\change{Applying our parametric dynamic model, we found a significant treatment effect with $p$- value
$p=1.0158\times10^{-5}$,
indicating that the customer subsidy increased overall driver income. Figure \ref{fig:realdata} presents the fitted values and relative residuals for two regions whose driver incomes dominate the overall city revenue. The relative residuals are computed by dividing the fitted residuals by the standard deviation of the corresponding regional incomes. The fitted patterns closely track the observed daily income dynamics, and the residuals show no systematic deviation, indicating a good model fit in economically active areas. In contrast, a few low-income regions with sparse ride requests exhibit larger residual fluctuations, likely due to weaker signal strength and higher relative noise. Given their limited contribution to the aggregate outcome, we focus on the dominant regions for clarity and space economy. This example demonstrates how our framework can evaluate dynamic spatial policies in large-scale online experiments.}

\section{Concluding Remarks}\label{sec:discussion}
In this paper, we aim at enhancing  the efficiency of inference across various scenarios by assigning treatment randomly across units. A pivotal discovery  is that randomization effectively diminishes the spatial correlation of noise. Under standard conditions, we observe that the MSE ratio of ATE estimators, derived from spatially randomized designs compared to classical global designs, is inversely related to the number of units or clusters, applicable in both static and dynamic frameworks. %Additionally, we advocate for the application of smoothing techniques on OLS coefficient estimators as a means to lower MSE, although this does not enhance the order of convergence of the outcomes and thus is omitted from the main discussion. Further details on this method are available in the supplementary materials.

Moving forward, there are several important topics worthy of further investigation. Firstly, examining
experimental designs with temporal random effects is an interesting issue, which may cause endogeneity bias in the state autoregression model.
% , resulting in inconsistent ATE estimation. 
% To tackle this inconsistency,
One %approach is to use historical data in which actions were set to baseline policies to estimate autoregressive coefficients. Another 
possible solution is to assume that random effects satisfy certain covariance structures such as declining correlation as temporal distance increases. 
Secondly, we have focused on the inference properties of multi-stages with finite time frames. For infinite
stages, where a variance estimate is likely to diverge, Gaussian approximation techniques
may be applied to overcome the problem, as discussed in \cite{luo2022policy}. 
% Thirdly, while \cite{lei2021regression} examined inference efficiency with high-dimensional state variables under a completely randomized design, the challenges posed by spatially randomized designs are worthy of investigation.
Finally, 
% we considered
% a general and simplistic interference structure in this study. I
investigating more complex interference structures and/or more complex treatment variables \citep{ao2021multivalued,dong2023regression} would be an interesting problem.

% \section*{APPENDICES MUST BE MOVED TO A SUPPLEMENT FILE}
%%%%%%%%%%%%%%%%%%%%%%%%%%%%%%%%%%%%%%%%%%%%%%
%% Support information, if any,             %%
%% should be provided in the                %%
%% Acknowledgements section.                %%
%%%%%%%%%%%%%%%%%%%%%%%%%%%%%%%%%%%%%%%%%%%%%%
\begin{acks}[Acknowledgments]
% This research is supported by National Key R\&D Program of China grants (No. 2022YFA1003800), the National Natural Science Foundation of China (No. 72301276, 12292981, 11931001), the LMAM and the Fundamental Research Funds for the Central Universities (LMEQF) and an EPSRC grant EP/W014971/1. 
Hongtu Zhu is the corresponding author. E-mail: htzhu@email.unc.edu.  
\end{acks}

%%%%%%%%%%%%%%%%%%%%%%%%%%%%%%%%%%%%%%%%%%%%%%
%% Funding information, if any,             %%
%% should be provided in the                %%
%% funding section.                         %%
%%%%%%%%%%%%%%%%%%%%%%%%%%%%%%%%%%%%%%%%%%%%%%
\begin{funding}
This research is supported by National Key R\&D Program of China grants (No. 2022YFA1003800), the National Natural Science Foundation of China (No. 72301276, 12571309, 62588101, 12292981, 11931001), the LMAM and the Fundamental Research Funds for the Central Universities (LMEQF) and an EPSRC grant EP/W014971/1.
\end{funding}

%%%%%%%%%%%%%%%%%%%%%%%%%%%%%%%%%%%%%%%%%%%%%%
%% Supplementary Material, including data   %%
%% sets and code, should be provided in     %%
%% {supplement} environment with title      %%
%% and short description. It cannot be      %%
%% available exclusively as external link.  %%
%% All Supplementary Material must be       %%
%% available to the reader on Project       %%
%% Euclid with the published article.       %%
%%%%%%%%%%%%%%%%%%%%%%%%%%%%%%%%%%%%%%%%%%%%%%
\begin{supplement}
The supplementary materials contain the additional implementation details, simulation results and supplementary theoretical results along with the proofs. The code and data can be found on the GitHub page at https://github.com/anneyang0060/SpatiallyRandomization. \change{The real data used in our study is proprietary and cannot be shared publicly. However, we have provided a synthetic dataset that can yield similar results for broader accessibility and understanding.}
\end{supplement}

%%  imsart-nameyear.bst  will be used to                   %%
%%  MR numbers will be added by VTeX.                      %%
%%                                                         %%
%%  Use \cite{...} to cite references in text.             %%
%%                                                         %%
%%%%%%%%%%%%%%%%%%%%%%%%%%%%%%%%%%%%%%%%%%%%%%%%%%%%%%%%%%%%%

%% if your bibliography is in bibtex format, uncomment commands:
%\bibliographystyle{imsart-nameyear} % Style BST file
%\bibliography{bibliography}       % Bibliography file (usually '*.bib')

%% or include bibliography directly:
\bibliographystyle{imsart-nameyear} % Style BST file (imsart-number.bst or imsart-nameyear.bst)
\bibliography{ref}       % Bibliography file (usually '*.bib')

\end{document}